\documentclass[a4paper]{article}[11pt]
\usepackage{amsfonts}
\usepackage{amsmath,cite,tikz,multirow,amssymb,mathrsfs,latexsym,bm,inputenc,longtable,amsfonts}
\usepackage{tabularx}
\usepackage{array}
\usepackage{float}
\usepackage[english]{babel}
\usepackage{enumitem}
\usepackage{colortbl}
\usepackage{xcolor}
\usepackage{hhline}
\begin{document}

\newcommand{\qed}{\hphantom{.}\hfill $\Box$\medbreak}
\newcommand{\Proof}{\noindent{\bf Proof \ }}

\newtheorem{theorem}{Theorem}[section]
\newtheorem{lemma}[theorem]{Lemma}
\newtheorem{corollary}[theorem]{Corollary}
\newtheorem{remark}[theorem]{Remark}
\newtheorem{example}[theorem]{Example}
\newtheorem{definition}[theorem]{Definition}
\newtheorem{construction}[theorem]{Construction}
\newtheorem{proposition}[theorem]{Proposition}

\title{\large{\bf On the cardinalities of quantum Latin squares  \footnote{Supported by the Natural Science Foundation of Hebei Province   Grant No. A2025205023, the National Natural Science Foundation of China  Grant No. 11871019 and Science Foundation of Hebei Normal University   Grant No. L2024B52.}}}

\author{Yajuan Zang$^{a}$, Meihui Zheng$^{b}$, Zihong Tian$^{b}$\thanks{Corresponding author. E-mail address: tianzh68@163.com}, Xiuling Shan$^{b}$ \\
\small $^a$ Department of Elementary Education, Hebei Normal University,
          Shijiazhuang 050024,  P. R. China\\
\small$^b$ School of Mathematical Sciences, Hebei Normal University,
          Shijiazhuang 050024, P. R. China }

\date{}

\maketitle

\noindent {\bf Abstract:} A quantum Latin square of order $v$, QLS($v$), is a $v\times v$ array in which each of entries is  a unit column vector from the Hilbert space $\mathbb{C}^{v}$, such that every row and column forms an orthonormal basis of $\mathbb{C}^{v}$. The  cardinality  of a QLS($v$) is the number of its vectors distinct up to a global phase, which is the  crucial indicator for distinguishing between   classical QLSs and non-classical QLSs. In this paper, we investigate the possible cardinalities of a QLS($v$).  As a result, we completely resolve the existence of a QLS($v$) with maximal cardinality for any $v\geq 4$. Moreover, based on Wilson's construction and Direct Product construction,  we establish some possible cardinality range  of a QLS($v$) for any $v\geq 4$. \vspace{0.1cm}


\section{Introduction}

Latin square is an important combinatorial configuration that constitute a highly active research, with signification applications in diverse fields such as statistics \cite{Gao}, cryptography \cite{Pal,Laywine}, computer science \cite{Donald}, communications \cite{Hayashi} and bioinformatics \cite{Ryazanov}.

Let $[v]=\{0,1,\ldots,v-1\}$. A \emph{Latin square} of order $v$, LS$(v)$, is a $v\times v$ array $L=(L_{i,j})$ in which each of the elements of $[v]$ occurs exactly once in each row and column. A Latin square $L$ of order $v$ is \emph{idempotent} if $L(i,i)=i$ for $i\in [v]$. Two Latin squares $L_1$, $L_2$ of order $v$ are \emph{orthogonal}, denoted by OLS$(v)$, if when $L_1$ is superimposed on $L_2$, every ordered pair $(0,0), (0,1), \dots, (v-1,v-1)$ occurs. It is well known that an LS($v$) exists for any $v\geq 2$, an idempotent LS($v$) exists for any $v\geq 3$ and a pair of OLS($v$)s exists  for any $v\geq 3$ and $v\neq 6$\cite{Colbourn}.

Quantum Latin squares are the quantum analogues of  Latin squares, introduced by Musto and Vicary \cite{Musto2016}. They are connected to various other primitives used in quantum information, such as unitary error bases  \cite{Musto2016}, mutually unbiased bases  \cite{Musto2017,Zang2022}, $k$-uniform states  \cite{Goyeneche2018,Zang2021}, controlled families of Hadamards   \cite{Reutter2019}, quantum teleportation and quantum error correction \cite{Helwig2012}, etc.

A \emph{quantum~Latin~square} of order $v$, QLS($v$), is a $v\times v$ array $\Phi=(|\Phi_{i,j}\rangle)$ in which each of entries is a unit column vector  from the Hilbert space $\mathbb{C}^{v}$, such that every row and column forms an orthonormal basis of $\mathbb{C}^{v}$.
Specifically, by replacing each number $i\in [v]$ in a Latin square of order $v$ with the unit vector $|i\rangle\in \mathbb{C}^v$, where $|i\rangle$ is a unit column vector with its $(i+1)$ component equal to 1 and the others equal to 0, a QLS($v$) can be obtained. We call the orthonormal basis $\{|0\rangle,|1\rangle,\ldots,|v-1\rangle\}$ of $\mathbb{C}^v$ a \emph{computational basis}.

Two quantum~Latin~squares  $\Phi$, $\Psi$ of order $v$ are \emph{equivalent} if there exists some unitary operator $U$ on~$\mathbb{C}^v$, a family of modulus 1 complex numbers~$c_{ij}$ and permutations $\sigma, \tau \in S_v$, such that
\begin{equation}
|\Psi_{i,j}\rangle=c_{ij}U|\Phi_{\sigma(i), \tau(j)}\rangle,
\end{equation} where $i,j \in [v]$, $S_v$ is the symmetry group of order $v$.

Two unit vectors $|u\rangle,|v\rangle \in \mathbb{C}^{v}$ are said having the \emph{same global phase}, if there exists a $\theta \in [0,2\pi)$ such that
\begin{equation}
|u\rangle=\mathrm{e}^{\mathrm{i}\theta}|v\rangle,
\end{equation}
where $\mathrm{i}$ is the imaginary unit. Otherwise, they are considered \emph{distinct}. By the Cauchy-Schwarz inequality, we know that two unit vectors $|u\rangle, |v\rangle \in \mathbb{C}^{v}$ have the same global phase if and only if $|\langle u | v\rangle|=1$. In contrast, if $|u\rangle$ and $|v\rangle$ have distinct global phases, then $|\langle u | v \rangle| \in [0, 1)$. The \emph{cardinality} $c$ of a QLS($v$) is the number of its vectors distinct up to a global phase. It is not difficult to see that the cardinality $c$ of a QLS($v$) satisfies $v\leq c\leq v^2$.

Obviously, when the QLS($v$) is directly from an LS($v$), i.e., it is consisted by computational basis, then $c=v$. Such QLS($v$) is called  \emph{classical}. Moreover, the QLS($v$) which is equivalent with a classical QLS($v$) is also called \emph{classical} (or \emph{apparently quantum}).  Otherwise, a QLS($v$) is called \emph{non-classical} (or \emph{genuinely quantum}). It is easy to see that the cardinality of a non-classical QLS($v$)  satisfies $c>v$. When $v=2,3$,  Paczos et al.  \cite{Paczos} proved that a non-classical QLS($v$) cannot exist. When $v\geq 4$, Han et al.  \cite{Han} showed that a non-classical QLS($v$) exists.

The cardinality is the crucial indicator for distinguishing between   classical QLSs and non-classical QLSs.   Paczos et al. \cite{Paczos} analyzed all possible cardinalities of a quantum Sudoku square of order 4, which is a special  QLS just like the relationship  of Sudoku square and LS,   are 4, 6, 8 and 16. Zang et al. \cite{Zang2021} introduced several construction methods to get non-classical QLSs and  non-classical mutually orthogonal QLSs. Recently, Cao  et al. \cite{Cao} established the existence of QLS($v$) with maximal cardinality for all orders $v\geq 9\times 77^4+4$. Ji  et al.\cite{Ji} constructed QLS($4m$)s with all possible cardinalities for $m\geq 2$.

In this paper, we investigate the possible cardinalities of a QLS($v$). As a result, we completely resolve the existence of a QLS($v$) with maximal cardinality for any $v\geq 4$. Moreover, we discuss the cardinality range  of a QLS($v$) for any $v\geq 4$.

The remainder of this paper is organized as follows. In Section 2, we introduce the concept of cardinality for two unitary matrices. In addition, we present the construction method for infinite   unitary matrices in which any two matrices have distinct cardinalities. In Section 3, by using Fourier matrix, we construct  QLS($v$)s with maximal cardinality for any $v\geq 4$. In Section 4, we establish Wilson's construction and Direct Product construction for QLS($v$)s. Consequently,  we investigate the possible cardinality range  of  a QLS($v$)  for any $v\geq 4$ and $v\neq 5$. In Section 5, we summarize the main results presented in this paper.

\section{Preliminaries}
In this section, we introduce the concept of cardinality of two unitary matrices for the needs of subsequent research. Furthermore, we establish the existence of  unitary matrices with distinct cardinalities.

Before proceeding, we give some notations. Let $F_v=(F_v(k,l))$ be Fourier matrix of order $v$ as follows,
\begin{eqnarray}\footnotesize
F_v=\frac{1}{\sqrt{v}}\left(\begin{array}{ccccc}
1& 1&1 & 1&1\\
1  & \omega &\omega^2 &\cdots &\omega^{v-1} \\
1  & \omega^2 &\omega^4 &\cdots &\omega^{2(v-1)} \\
1  & \vdots &\vdots &\ddots &\vdots\\
1  & \omega^{v-1} &\omega^{2(v-1)}&\cdots &\omega^{(v-1)^2} \\
\end{array}\right),\label{e0}
\end{eqnarray}
where $F_v(k,l)=\frac{1}{\sqrt v} \omega^{kl}$, $\omega=\mathrm{e}^{\frac{2\pi \mathrm{i}}{v}}$, $k,l \in [v]$. $I_v$ is the identity matrix of order $v$. For positive integers $a, b$, $[a,b]=\{a, a+1,\dots, b\}$, $[a,b]_o$ and $[a,b]_e$ are the sets of odd, even numbers in $[a,b]$ respectively.  If $i\neq j$, $\delta_{ij}=0$, otherwise, $\delta_{ij}=1$.

Two unitary matrices $U$, $V$ of order $v$ are said having the \emph{same global phases}, if for any column vector in $U$, there just exists one column vector in $V$ with the same global phase.
A unitary matrix $U$ is said having \emph{distinct global phases} with unitary matrix $V$, if there exists at least  one  column vector in $U$ that has different global phase with all column vectors in $V$.
Two unitary matrices $U$, $V$ of order $v$ are said having  \emph{totally distinct  global phases}, if any column vector in $U$ and any column vector in $V$ have distinct global phases.

For unitary matrices $U$, $V$ of order $v$, the \emph{cardinality} $C$ of $U$ with $V$ is the number of the column vectors in $U$ distinct with all the column vectors in $V$ up to global phase.

\begin{lemma}\label{u0}
For any unitary matrices $U$, $V$ of order $v$, $0 \leq C\leq v$ and $C\neq 1$.
\end{lemma}

\Proof Here we just prove the latter, since the former is  clearly valid. Suppose $U=(\alpha_0~\alpha_1~\cdots~\alpha_{v-1})$ and $V=(\beta_0~\beta_1~\cdots~\beta_{v-1})$, where $\alpha_i$ and $\beta_j$ are the $(i+1)$-th, $(j+1)$-th column vectors of $U, V$ respectively, $i,j\in[v]$. Without loss of generality, assume $\alpha_{v-1}$ has  different global phase with all column vectors in $V$, and for any $i \in [v] \backslash \{v-1\}$, $\alpha_i$ and $\beta_i$ have the same global phase. Let $W_1 = span\{{{\alpha}_{0}}, {{\alpha}_{1}}, \dots , {{\alpha}_{v-2}}\} =span\{{{\beta}_{0}}, {{\beta}_{1}}, \dots , {{\beta}_{v-2}}\}$ and $\mathbb{C}^{v} = W_1 \bigoplus {W_2}$,
where $W_2$ is the orthogonal complement space of $W_1$. Therefore $W_{2} = span\{\alpha_{v-1}\}
=  span\{\beta_{v-1}\},$ which is a contradiction. \qed

According to the above description,  the unitary matrices $U$ and $V$ have the same global phase when $C=0$;  the unitary matrix $U$ has distinct global phases with unitary matrix $V$ when $C>0$; the unitary matrices $U$ and $V$ have totally distinct global phases when $C=v$.

In the following, we discuss the distinct cardinalities for some special unitary matrices. Let
\begin{equation}\footnotesize
E^\theta=\left(\begin{array}{cccc}
\mathrm{e}^{\mathrm{i}\theta}& & & \\
   & 1 & & \\
  &  &\ddots &\\
 & & & 1 \\
\end{array}\right)\label{e1}
\end{equation}
be a $v\times v$ matrix with $\theta\in [0,2\pi)$. Define
\begin{equation}
F_v^\theta=E^\theta F_v, \label{f1}
\end{equation}
where $F_v$ is the Fourier matrix of order $v$.  It is obvious that $F_v^\theta$ is a unitary matrix.

\begin{lemma}\label{u1}
For any positive integer $v \geq 2$, there exist infinite unitary matrices of order $v$ with totally distinct global phases.
\end{lemma}

\Proof When $\theta_1 \neq \theta_2$, $\theta_1,\theta_2\in [0,\pi)$,  we have $F_v^{\theta_1}$ and $F_v^{\theta_2}$ defined in Eq.(\ref{f1}) have totally distinct  global phases. In fact, suppose $F_v^{\theta_1}=(\alpha_0~ \alpha_1~ \cdots ~\alpha_{v-1})$, $F_v^{\theta_2}=(\beta_0~ \beta_1~ \cdots~ \beta_{v-1})$, where $\alpha_i$ and $\beta_j$ are the $(i+1)$-th, $(j+1)$-th column vectors of $U, V$ respectively, $i,j\in[v]$. For any $i,j \in [v]$,
\begin{eqnarray*}\footnotesize
&(\alpha_i, \beta_j)=\frac{1}{v}(e^{-\mathrm{i}\theta_1}~\omega^{-i}~\omega^{-2i}~\cdots~ \omega^{-(v-1)i})\left(\begin{array}{c}
 e^{\mathrm{i}\theta_2} \\
 \omega^{j}\\
  \omega^{2j}\\
  \vdots\\
  \omega^{(v-1)j}\\
\end{array}\right)\\
&\hspace{-1.7cm}=\frac{1}{v}(e^{\mathrm{i}(\theta_2-\theta_1)} + \sum\limits_{k=1}^{v-1} \omega ^{k(j-i)}).
\end{eqnarray*}

\noindent When $i\neq j$, since $\sum_{k=0}^{v-1} \omega ^{k(j-i)} = 0$, so  $|(\alpha_i, \beta_j)|=|\frac{1}{v}(e^{\mathrm{i}(\theta_2-\theta_1)}-1)| \neq 1$;
when $i= j$, $|(\alpha_i, \beta_j)|=|\frac{1}{v}(e^{i(\theta_2-\theta_1)}+(v-1))| \neq 1$. Consequently, if $\theta_1\neq \theta_2$,  $F^{\theta_1}_v$ and $F^{\theta_2}_v$ have totally distinct  global phases.

Let
\begin{equation}
\mathcal{F}=\{F^{\theta_i}_v |~ \theta_i\in [0,\pi), i\in \mathbb{N}^{+} \}
\end{equation}
be a family of unitary matrices. For any $F^{\theta_i}_v , F^{\theta_j}_v \in \mathcal{F}$, when $\theta_i\neq \theta_j$, then $F^{\theta_i}_v$ and $F^{\theta_j}_v$ have totally distinct  global phases, which ends the proof.  \qed

Next, we define a new  unitary matrix of order $v$ as follows,
\begin{eqnarray}\footnotesize
U=\begin{pmatrix}
1& & &\\
&  \ddots & & \\
&  & 1 & \\
& & &F^\theta_s
\end{pmatrix},\label{f2}
\end{eqnarray}
where $F^\theta_s$ is defined in Eq.(\ref{f1}), $2\leq s\leq v$, $\theta\in [0,\pi)$. Obviously, the matrix $U$ is unitary. Moreover, it is easy to see that for any column vector from the last $s$ columns in $U$ has distinct global phases with all the column vectors of the identity matrix  $I_v$. Thus we have the following conclusion.

\begin{lemma}\label{u2}
For any positive integer $v \geq 2$,  $2\leq s\leq v$,  there exists a unitary matrix $U$ of order $v$ such that the cardinality   of $U$ with $I_v$ is $s$.
\end{lemma}

\begin{lemma}\label{u3}
Suppose the matrix  $U$ of order $v$ is defined in $\rm{Eq.(\ref{f2})}$.  For any positive integers $v \geq 2$, $s\leq t\leq v$, there exists a unitary matrix $V$ of order $v$ such that the cardinality   of $V$ with $U$ is $t$.
\end{lemma}

\Proof Let
\begin{eqnarray}\footnotesize
V=\begin{pmatrix}
1& & &\\
&  \ddots & & \\
&  & 1 & \\
& & &F^{\theta'}_t
\end{pmatrix}
\end{eqnarray}
be a $v\times v $ matrix, where $F^{\theta'}_t$ is defined in Eq.(\ref{f1}), $\theta\neq\theta'$, $s\leq t\leq v$, $\theta,\theta'\in [0,\pi)$.

When $t=s$, on the one hand, arbitrary column vector from the last $t$ columns in $V$ is orthogonal with any vector from the first $v-t$ columns in $U$. On the other hand, arbitrary column vector from the last $t$ columns in $V$ has distinct global phases with any vector from the last $t$ columns in $U$ via Lemma \ref{u1}. In the same way, it is easy to prove  the cases of  $t>s$.
\qed

From the Lemma \ref{u1}-Lemma \ref{u3}, we have the corollary below.
\begin{corollary}\label{u4}
For any positive integers $v \geq 2$, $l\geq 1$, and the given $ l$ unitary matrices
$$\footnotesize U_i=\begin{pmatrix}
1& & &\\
&  \ddots & & \\
&  & 1 & \\
& & &F^{\theta_i}_{s_i}
\end{pmatrix},$$
where $F^{\theta_i}_{s_i}$ is defined in $\rm{Eq.(\ref{f1})}$, $2\leq s_i\leq v$,  $\theta_i \in [0,\pi)$, $1\leq i\leq l$,  there exists a unitary matrix  $$\footnotesize V=\begin{pmatrix}
1& & &\\
&  \ddots & & \\
&  & 1 & \\
& & &F^{\theta'}_t
\end{pmatrix},$$
where $\theta'\neq \theta_i$, $1\leq i\leq l$,  $\theta'\in [0,\pi)$, such that for any $2\leq t\leq v$ and $1\leq i\leq l$, the cardinality of $V$ with $U_i$ and $I_v$ is $t$.
\end{corollary}

\section{Quantum Latin squares with maximal cardinality}

In this section, we utilize Fourier matrix to construct QLS$(v)$s with maximal cardinality for any $v\geq 4$ .

\begin{theorem}\label{cc}
For any positive integer $v\geq4$, there exists a $\rm{QLS}$$(v)$ with~maximal~cardinality.
\end{theorem}

\Proof For any positive integer $v\geq4$, let $F_v$ be the Fourier matrix of order $v$. Assume $\hat{F}_v = \sqrt{v}{{F}_{v}}$. Define a new unitary matrix $M$ of order $v$ by exchanging the second and the third row of $\hat{F}_v $ as below
\begin{eqnarray}\footnotesize
M=\left(\begin{array}{ccccc}
1& 1&1 & 1&1\\
1  & \omega^2 &\omega^4 &\cdots &\omega^{2(v-1)} \\
1  & \omega &\omega^2 &\cdots &\omega^{v-1} \\
1  & \omega^3 &\omega^6 &\cdots &\omega^{3(v-1)} \\
1  & \vdots &\vdots &\ddots &\vdots\\
1  & \omega^{v-1} &\omega^{2(v-1)}&\cdots &\omega^{(v-1)^2} \\
\end{array}\right)=(M_0~M_1~\cdots~M_{v-1}),\label{m}
\end{eqnarray}
where $M_i$ is the $(i+1)$-th column vector of $M$. For any $i\in[v]$, assume
\begin{eqnarray}\footnotesize
M^*_{i}=\begin{pmatrix}
1& & & & & \\
& \omega^{2i} & & & &\\
& & \omega^{i} &  & & \\
& & &\omega^{3i} & & \\
& & & & \ddots & \\
& & & &   &\omega^{(v-1)i}
\end{pmatrix}.\label{u11}
\end{eqnarray}
Define a new matrix
 \begin{equation}
U_i= M^*_i \cdot F_v \label{uu}
\end{equation}
 with
 $$U_i=(\alpha^{i}_{0}~\alpha^{i}_{1}~\cdots~\alpha^{i}_{v-1}),$$
 where $\alpha^{i}_{j}$ is the $(j+1)$-th column vector of $U_i$, $j\in[v]$. It is easy to see that for any $i\in [v]$, $U_i$ is a unitary matrix and $U_0=F_v$.

Define $\Phi=(|\Phi_{i,j}\rangle)_{i,j \in [v]}$ with
$|\Phi_{i,j}\rangle = \alpha^{i}_{j}$, $i,j\in[v]$. Then we claim that $\Phi$ is a QLS($v$) with $c=v^2$.

 First, we prove that $\Phi$ is a QLS($v$).
Fixed $i\in[v]$, for any $j,j'\in[v]$, since $U_{i}$ is a unitary matrix, thus $ \langle \Phi_{i,j}|\Phi_{i,j'} \rangle=(\alpha_j^i,\alpha_{j'}^i)=\delta_{jj'}$.
Fixed $j\in[v]$, for any $i,i'\in[v]$, since

\hspace{1.5cm}$\footnotesize \centering\small {\begin{tabular}{cc}
$|\Phi_{i,j}\rangle=\alpha^{i}_{j}=M^*_i \cdot \alpha^{0}_{j}= \frac{1}{\sqrt{v}}
\begin{pmatrix} 1 \\ \omega^{2i+j} \\ \omega^{i+2j} \\ \omega^{3(i+j)} \\ \vdots \\ \omega^{(v-1)(i+j)} \end{pmatrix}$,~~
$|\Phi_{i',j}\rangle= \frac{1}{\sqrt{v}} \begin{pmatrix} 1 \\ \omega^{2i'+j} \\ \omega^{i'+2j} \\ \omega^{3(i'+j)} \\ \vdots \\ \omega^{(v-1)(i'+j)} \end{pmatrix}$,
\end{tabular}}$

\noindent we have $\langle \Phi_{i,j}|\Phi_{i',j} \rangle=\frac{1}{v}\sum\limits_{k=0}^{v-1}\omega^{k(i'-i)}=\delta_{ii'}$.

Next, we show that the cardinality of $\Phi$ is $v^{2}$.
For any $i,i',j,j' \in [v], (i,j)\neq(i',j')$,  we consider the following two cases.

When  $i=i'$ or $j=j'$, we have $\langle \Phi_{i,j}|\Phi_{i',j'}\rangle = 0$ from above. Thus $|\Phi_{i,j}\rangle$ and $|\Phi_{i',j'}\rangle$ have distinct global phases.

When $i \neq i'$ and $ j \neq j'$, suppose that $|\Phi_{i,j}\rangle$, $|\Phi_{i',j'}\rangle$ have the same global phase, then the equation $|\Phi_{i,j}\rangle$ = $|\Phi_{i',j'}\rangle$ is established by the first component of them. Thus we have
\begin{equation*}
\left\{
\begin{array}{c}
\omega^{i+2j} = \omega^{i'+2j'};\\
\omega^{(v-1)(i+j)} = \omega^{(v-1)(i'+j')},
\end{array}
\right.
\end{equation*}
i.e.,
\begin{equation*}
\left\{
\begin{array}{c}
i+2j\equiv i'+2j' ~(\text{mod}~v);\\
i+j \equiv i'+j' ~(\text{mod}~v).
\end{array}
\right.
\end{equation*}
Therefore,  $i \equiv i' (\text{mod}~v), j \equiv j' (\text{mod}~ v)$, which is a contradiction. \qed

\begin{example}
There exists a $\rm{QLS(5)}$ with maximal cardinality $25$.
\end{example}
\Proof Let $\hat{F}_5=\sqrt{5}F_5$ and $M$ be the unitary matrix by exchanging the second and the third row of $\hat{F}_5$ as below
\begin{equation}\footnotesize
M =
\small{\begin{pmatrix}
1 & 1 & 1 & 1 & 1\\
1 & \omega^{2} & \omega^{4} & \omega & \omega^{3}\\
1 & \omega & \omega^{2} & \omega^{3} & \omega^{4}\\
1 & \omega^{3} & \omega & \omega^{4} & \omega^{2}\\
1 & \omega^{4} & \omega^{3} & \omega^{2} & \omega
\end{pmatrix}.}
\end{equation}
For any $i\in [5]$, suppose $M^*_i$ is the diagonal unitary matrix defined from the $(i+1)$-th column vector of $M$ as Eq.(\ref{u11}). Define
\begin{equation*}\footnotesize
U_{0} = M^*_0 \cdot F_5= F_{5} =\frac{1}{\sqrt5}
\begin{pmatrix}
1 & 1 & 1 & 1 & 1\\
1 & \omega & \omega^{2} & \omega^{3} & \omega^{4}\\
1 & \omega^{2} & \omega^{4} & \omega & \omega^{3}\\
1 & \omega^{3} & \omega & \omega^{4} & \omega^{2}\\
1 & \omega^{4} & \omega^{3} & \omega^{2} & \omega
\end{pmatrix};
\end{equation*}
\begin{equation*}
\centering\scalebox{0.9} {\begin{tabular}{cc}
$U_{1} = M^*_1 \cdot F_5=
\frac{1}{\sqrt{5}}
\begin{pmatrix}
1 & 1 & 1 & 1 & 1\\
\omega^{2} & \omega^{3} & \omega^{4} & 1 & \omega\\
\omega & \omega^{3} & 1 & \omega^{2} & \omega^{4}\\
\omega^{3} & \omega & \omega^{4} & \omega^{2} & 1\\
\omega^{4} & \omega^{3} & \omega^{2} & \omega & 1
\end{pmatrix};~
U_{2} = M^*_2 \cdot F_5  =
\frac{1}{\sqrt{5}}
\begin{pmatrix}
1 & 1 & 1 & 1 & 1\\
\omega^{4} & 1 & \omega & \omega^{2} & \omega^{3}\\
\omega^{2} & \omega^{4} & \omega & \omega^{3} & 1\\
\omega & \omega^{4} & \omega^{2} & 1 & \omega^{3}\\
\omega^{3} & \omega^{2} & \omega & 1 & \omega^{4}
\end{pmatrix};$
\end{tabular}}
\end{equation*}
\begin{equation*}
\centering\scalebox{0.9} {\begin{tabular}{cc}
$U_{3} = M^*_3 \cdot F_5 =
\frac{1}{\sqrt{5}}
\begin{pmatrix}
1 & 1 & 1 & 1 & 1\\
\omega & \omega^{2} & \omega^{3} & \omega^{4} & 1\\
\omega^{3} & 1 & \omega^{2} & \omega^{4} & \omega\\
\omega^{4} & \omega^{2} & 1 & \omega^{3} & \omega\\
\omega^{2} & \omega & 1 & \omega^{4} & \omega^{3}
\end{pmatrix};~
U_{4} = M^*_4 \cdot F_5 =
\frac{1}{\sqrt{5}}
\begin{pmatrix}
1 & 1 & 1 & 1 & 1\\
\omega^{3} & \omega^{4} & 1 & \omega & \omega^{2}\\
\omega^{4} & \omega & \omega^{3} & 1 & \omega^{2}\\
\omega^{2} & 1 & \omega^{3} & \omega & \omega^{4}\\
\omega & 1 & \omega^{4} & \omega^{3} & \omega^{2}
\end{pmatrix}.$
 \end{tabular}}
\end{equation*}
\noindent Consequently, a QLS(5)  with maximal cardinality is given.
\begin{equation*}
\centering\scalebox{0.69} {\begin{tabular}{c}
$\begin{array}{|c|c|c|c|c|}
  \hline
\frac{|0 \rangle+|1 \rangle+|2 \rangle+|3 \rangle+|4 \rangle}{\sqrt{5}}
& \frac{|0 \rangle+\omega|1 \rangle+\omega^2|2 \rangle+\omega^3|3 \rangle+\omega^4|4 \rangle}{\sqrt{5}}
& \frac{|0 \rangle+\omega^2|1 \rangle+\omega^4|2 \rangle+\omega|3 \rangle+\omega^3|4 \rangle}{\sqrt{5}}
& \frac{|0 \rangle+\omega^3|1 \rangle+\omega|2 \rangle+\omega^4|3 \rangle+\omega^2|4 \rangle}{\sqrt{5}}
& \frac{|0 \rangle+\omega^4|1 \rangle+\omega^3|2 \rangle+\omega^2|3 \rangle+\omega|4 \rangle}{\sqrt{5}} \\
  \hline
\frac{|0 \rangle+\omega^2|1 \rangle+\omega|2 \rangle+\omega^3|3 \rangle+\omega^4|4 \rangle}{\sqrt{5}}
& \frac{|0 \rangle+\omega^3|1 \rangle+\omega^3|2 \rangle+\omega|3 \rangle+\omega^3|4 \rangle}{\sqrt{5}}
& \frac{|0 \rangle+\omega^4|1 \rangle+|2 \rangle+\omega^4|3 \rangle+\omega^2|4 \rangle}{\sqrt{5}}
& \frac{|0 \rangle+|1 \rangle+\omega^2|2 \rangle+\omega^2|3 \rangle+\omega|4 \rangle}{\sqrt{5}}
& \frac{|0 \rangle+\omega|1 \rangle+\omega^4|2 \rangle+|3 \rangle+|4 \rangle}{\sqrt{5}} \\
  \hline
\frac{|0 \rangle+\omega^4|1 \rangle+\omega^2|2 \rangle+\omega|3 \rangle+\omega^3|4 \rangle}{\sqrt{5}}
& \frac{|0 \rangle+|1 \rangle+\omega^4|2 \rangle+\omega^4|3 \rangle+\omega^2|4 \rangle}{\sqrt{5}}
& \frac{|0 \rangle+\omega|1 \rangle+\omega|2 \rangle+\omega^2|3 \rangle+\omega|4 \rangle}{\sqrt{5}}
& \frac{|0 \rangle+\omega^2|1 \rangle+\omega^3|2 \rangle+|3 \rangle+|4 \rangle}{\sqrt{5}}
& \frac{|0 \rangle+\omega^3|1 \rangle+|2 \rangle+\omega^3|3 \rangle+\omega^4|4 \rangle}{\sqrt{5}} \\
  \hline
\frac{|0 \rangle+\omega|1 \rangle+\omega^3|2 \rangle+\omega^4|3 \rangle+\omega^2|4 \rangle}{\sqrt{5}}
& \frac{|0 \rangle+\omega^2|1 \rangle+|2 \rangle+\omega^2|3 \rangle+\omega|4 \rangle}{\sqrt{5}}
& \frac{|0 \rangle+\omega^3|1 \rangle+\omega^2|2 \rangle+|3 \rangle+|4 \rangle}{\sqrt{5}}
& \frac{|0 \rangle+\omega^4|1 \rangle+\omega^4|2 \rangle+\omega^3|3 \rangle+\omega^4|4 \rangle}{\sqrt{5}}
& \frac{|0 \rangle+|1 \rangle+\omega|2 \rangle+\omega|3 \rangle+\omega^3|4 \rangle}{\sqrt{5}} \\
  \hline
\frac{|0 \rangle+\omega^3|1 \rangle+\omega^4|2 \rangle+\omega^2|3 \rangle+\omega|4 \rangle}{\sqrt{5}}
& \frac{|0 \rangle+\omega^4|1 \rangle+\omega|2 \rangle+|3 \rangle+|4 \rangle}{\sqrt{5}}
& \frac{|0 \rangle+|1 \rangle+\omega^3|2 \rangle+\omega^3|3 \rangle+\omega^4|4 \rangle}{\sqrt{5}}
& \frac{|0 \rangle+\omega|1 \rangle+|2 \rangle+\omega|3 \rangle+\omega^3|4 \rangle}{\sqrt{5}}
& \frac{|0 \rangle+\omega^2|1 \rangle+\omega^2|2 \rangle+\omega^4|3 \rangle+\omega^2|4 \rangle}{\sqrt{5}} \\
  \hline
\end{array}$
\end{tabular}}
\end{equation*}
\qed

\section{The possible cardinalities of quantum Latin squares}

In this section, we investigate the possible cardinality range  of  QLS$(v)$s. In combinatorial design theory, there are several  classical  and powerful recursive construction methods for Latin squares such as  Wilson's  construction, Direct Product construction, etc.  When $v=mt+s$, $s\geq 1$, we employ unitary matrices to implement Wilson's construction for quantum Latin squares and analyse the  cardinalities of QLS$(v)$s. When $v=mt$, we examine the cardinalities of QLS$(v)$s arising from the Direct Product construction described in the Ref. \cite{Zang2021}.

\begin{lemma}\label{c}
For any positive integer $v \geq 4$, the cardinality of a $\rm{QLS}$$(v)$ satisfies $c \neq v+1$.
\end{lemma}

\Proof Suppose $\Phi$ is a QLS($v$). Without loss of generality, we assume that the elements in the first row of $\Phi$ are from the computational basis of $\mathbb{C}^v$. If $c=v+1$, then only one element in $\Phi$ has the distinct global phases with all of the vectors in the computational basis. Suppose the element denoted by $\color{red}{\star}$ locates in the second row, first column of $\Phi$ as below.

\[
\begin{array}{|c|c|c|c|}
  \hline
  |0\rangle & |1\rangle & \dots & |v-1\rangle \\
  \hline
  \color{red}{\star}     &           &        &       \\
  \hline
  \vdots &   &   &   \\
  \hline
     &   &   &   \\
  \hline
\end{array}
\]

Define two matrices $U$ and $V$ which are constituted by the elements of the first row and second row in $\Phi$ as the column vectors respectively. It is easy to check that  $U$ and $V$ are unitary matrices.
From Lemma \ref{u0}, we know $U$ and $V$ have the same global phase, which is a contradiction. Thus the conclusion is established. \qed

\subsection{The cardinalities of QLS($mt+s$)s}
In this subsection, we generalize the Wilson's construction in Ref. \cite{Denes} to  generate quantum Latin squares and calculate their cardinalities.

\begin{construction}\label{qls1}
$\rm{(Wilson's~construction~for~QLS)}$ If there exists a pair of $\rm{OLS}$$(t)$s, an $\rm{LS}$$(s)$, an $\rm{LS}$$(m)$
and an $\rm{LS}$$(m+1)$, then there exists a $\rm{QLS}$$(mt+s)$, where $2\leq s<t$.
\end{construction}
\Proof
Let $(L(a,b))_{a,b\in [t]}$ and $(L'(a,b))_{a,b\in [t]}$ be a pair of OLS($t$)s, and define
\begin{eqnarray}
\mathcal{A}=\{(a,b,L(a,b),L'(a,b))~|~a,b\in [t]\}.
\end{eqnarray}

Let $\mathcal{B}=\{|0\rangle, |1\rangle,\ldots,|mt+s-1\rangle\}$ be the computational basis of the space $\mathbb{C}^{mt+s}$. Assume $R_a=[ma,ma+m-1]$ for any $a\in [t]$ and $S=[mt,mt+s-1]$.
Then $\big(\bigcup\limits_{a\in [t]}R_a\big) \bigcup S=[mt+s]$.
For any $a\in[t]$, assume
\begin{eqnarray}
W_a=span\{|i\rangle: i \in R_a\},\\
W_t=span\{|i\rangle: i \in S\}. ~~
\end{eqnarray}
Then $\big(\bigoplus\limits_{a\in [t]}W_a\big) \bigoplus W_t=\mathbb{C}^{mt+s}$.

Define some unitary matrices of order $mt+s$ as below
\begin{eqnarray}\footnotesize
\centering\small {\begin{tabular}{cc}\scalebox{1}{
$U_a =\begin{pmatrix}
U^0_a& & &\\
&  \ddots & & \\
&  & U^{t-1}_a & \\
& & &V_s
\end{pmatrix}, $~~ $U_t=
\begin{pmatrix}
I_m& & &\\
&  \ddots & & \\
&  & I_m & \\
& & &V'_s
\end{pmatrix}$},
 \end{tabular}}\label{f3}
 \end{eqnarray}
where  $U^i_a$ is a unitary matrix of order $m$, $a, i\in [t]$, $V_{s}$ and $V'_{s}$ are unitary matrices of order $s$.

Next, we construct an array  $\Phi=(|\Phi_{x,y}\rangle)$ of order $mt+s$.
For any $A \in \mathcal{A}$, $A=(a,b,L(a,b),\\L'(a,b))$, $a,b\in[t]$, we consider the entries of $\Phi$ based on  $A$ as follows.

(1) If $L'(a,b)\geq s$, since there is an LS($m$), then we can construct a QLS($m$), $(|\phi_{i,j}\rangle)_{i,j\in[m]}$, based on the orthonormal basis $\{U_a|mL(a,b)\rangle,U_a|mL(a,b)+1\rangle,\ldots,U_a|mL(a,b)+m-1\rangle\}$ of the subspace $W_{L(a,b)}$. For any $i,j\in[m]$, define
\begin{eqnarray}
|\Phi_{ma+i,mb+j}\rangle=|\phi_{i,j}\rangle.
\end{eqnarray}

(2) If $L'(a,b)< s$, since there is an LS($m+1$), without loss of generality, let the element $m$ appear  in the lower right corner. Then we can construct a QLS($m+1$), $(|\varphi_{i,j}\rangle)_{i,j\in[m+1]}$  with $|\varphi_{m,m}\rangle=|mt+L'(a,b)\rangle$ based on the orthonormal basis $\{U_a|mL(a,b)\rangle,\ldots,U_a|mL(a,b)+m-1\rangle, |mt+L'(a,b)\rangle\}$ of the subspace $W_{L(a,b)}\bigoplus span\{|mt+L'(a,b)\rangle\}$. For any $i,j\in[m]$, define
\begin{eqnarray}
|\Phi_{ma+i,mb+j}\rangle=|\varphi_{i,j}\rangle; ~~~~~~~\nonumber\\
|\Phi_{ma+i,mt+L'(a,b)}\rangle=|\varphi_{i,m}\rangle; \nonumber \\
|\Phi_{mt+L'(a,b),mb+j}\rangle=|\varphi_{m,j}\rangle.
\end{eqnarray}

(3) Since there is an LS($s$),   we can construct a QLS($s$), $(|\psi_{i,j}\rangle)_{i,j\in[s]}$, based on the orthonormal basis $\{U_t|mt\rangle, U_t|mt+1\rangle,\ldots, U_t|mt+s-1\rangle\}$ of the subspace $W_t$. For any $i,j\in [s]$, define
\begin{eqnarray}
|\Phi_{mt+i,mt+j}\rangle=|\psi_{i,j}\rangle.
\end{eqnarray}

We claim that $\Phi=(|\Phi_{x,y}\rangle)_{x,y\in[mt+s]}$ is a QLS($mt+s$).
In fact, for any $x,y\neq y'\in[mt+s]$, assuming $x=mg+i$, $y=mh+j$ and $y'=mh'+j'$, where $0\leq g,h,h'\leq t, 0\leq i,j,j'\leq max\{m,s-1\}$, it is necessary to prove the following equation is true,
\begin{eqnarray}\label{22}
\langle \Phi_{x,y}|\Phi_{x,y'}\rangle=0.
\end{eqnarray}

\textcircled{1} For $g\in[t]$, by the cases (1), (2) above, we need to consider the following two cases.

($i$) When $h=h'$, $j\neq j'$, if $h\in [t]$, then $|\Phi_{x,y}\rangle, |\Phi_{x,y'}\rangle$ are from the same row of $\phi$ or $\varphi$. So the Eq.(\ref{22}) is true. If $h=t$, there exist $b,b'\in[t]$ such that  $j=L'(g,b), j'=L'(g,b')$. Since $j\neq j'$, thus $b\neq b'$. So we have $W_{L(g,b)}\cap W_{L(g,b')}=\{0\}$, the Eq.(\ref{22}) is true.

($ii$) When $h\neq h'$, if $h,h'\in [t]$, since $L(g,h)\neq L(g,h')$ and $W_{L(g,h)}\cap W_{L(g,h')}=\{0\}$, so the Eq.(\ref{22}) is established. When one of $h,h'$ is $t$, without loss of generality, suppose $h=t$ and there is $b\in[t]$ such that $j=L'(g,b)$. Let $y=mt+L'(g,b)$, $y'=mh'+j', h'\in [t]$. If $h'=b$, then $|\Phi_{x,y}\rangle $ and $|\Phi_{x,y'}\rangle$ are from the same row of $\varphi$ based on the same subspace $W_{L(g,b)}$ or one is from the subspace $W_{L(g,b)}$ and the other is from the subspace $W_t$. Therefore, the Eq.(\ref{22}) is established. If $h'\neq b$, then $|\Phi_{x,y}\rangle$ and $ |\Phi_{x,y'}\rangle$ are from the subspaces $W_{L(g,b)}$ and $W_{L(g,h')}$ respectively or one is from $W_t$ and the other is from $W_{L(g,h')}$. The Eq.(\ref{22}) is true for the cases.

\textcircled{2} For $g=t$, by the cases (2), (3), we need to consider the following cases.

($i$) When $h,h'\neq t$, there are several ordered pairs $(a,b)$ satisfying $i=L'(a,b)$. Here we assume $i=L'(a,h)=L'(a',h')$ for $(a,h)\neq (a',h')$. Since $(L(a,b), L'(a,b))$ runs over $\mathbb{Z}_t\otimes \mathbb{Z}_t$ for $a,b\in [t]$, thus for fixed $i$,  $L(a,h)\neq L(a',h')$. So we have $|\Phi_{x,y}\rangle$ and $|\Phi_{x,y'}\rangle$ are from the same row of $\varphi$ based on the same subspace $W_{L(a,h)}$  or one is from the subspace $W_{L(a,h)}$ and the other is from the subspace $W_{L(a',h')}$. The Eq.(\ref{22}) is true.

($ii$) When one of $h,h'$ is $t$, then for $|\Phi_{x,y}\rangle$ and $|\Phi_{x,y'}\rangle$, one is from the subspace $W_{L(a,b)}$ and the other is from the subspace $W_t$. It is not difficult to see Eq.(\ref{22}) is true.

($iii$) When $h=h'=t$, then $|\Phi_{x,y}\rangle $ and $ |\Phi_{x,y'}\rangle$ are from the same row of $\psi$. Obviously, Eq.(\ref{22}) is true.

For any $y, x\neq x'\in[mt+s]$, assuming $y=mh+j$, $x=mg+i$ and $x'=mg'+i'$, where $0\leq g,g',h\leq t, 0\leq i,i',j\leq max\{m,s-1\}$, by the same way we can prove
\begin{eqnarray}
\langle \Phi_{x,y}|\Phi_{x',y}\rangle=0.
\end{eqnarray}

Thus   $\Phi$ is a QLS($mt+s$). \qed

In the following, we discuss the possible cardinalities of $\Phi$ constructed from Construction \ref{qls1}.

\begin{theorem}\label{c1}
For any $t\geq 2$, $t\neq 2, 6$ and $2 \leq s<t$, when $m\geq 2$, there exists a $\rm{QLS}$$(mt+s)$ with cardinality $c$ for any $c\in [mt+s,mt^2+2s]\setminus \{mt+s+1\}$, except possible for $c\in [2t+2,2t^2+4]_o$ with $m=2$ and $s=2$.
\end{theorem}

\Proof We consider the cardinalities of the quantum Latin square $\Phi$ constructed in Construction \ref{qls1}. Suppose $U_a$, $U_t$ are unitary matrices defined in Eq.(\ref{f3}), $a\in[t]$. Specially, we assume $U_0^i=I_m$, $i\in[t]$. For $1\leq a\leq t-1$, let $V_s=I_s$ and
\begin{eqnarray*}
 V_s'=\begin{pmatrix}
I_{s-j}& \\
& F_j^\theta
\end{pmatrix},~~ U_a^i=\begin{pmatrix}
I_{m-d_i}& \\
& F_{d_i}^{\theta_a}
\end{pmatrix},
\end{eqnarray*}
where $F_{d_i}^{\theta_a}$ and $F_j^\theta$  are defined in Eq.(\ref{f1}) for $2\leq d_i\leq m$, $\theta_a\neq \theta_{a'}\neq \theta \in[0, \pi)$, $1\leq a\neq a'\leq t-1$, $2\leq j\leq s$, $0\leq i\leq t-1$.

For $1\leq a\leq t-1$, suppose $U_a$ have $x_a$ column vectors which have totally distinct  global phases with $U_j$ for any $0\leq j \leq a-1$. Let $x'$ be the cardinality of $V_s'$ with $V_s$. Thus the cardinality of $\Phi$ is
\begin{equation}
c=\sum\limits_{a=1}^{t-1} x_a+x'+(mt+s).\label{e2}
\end{equation}
For $i\in[t]$, suppose $U_a^i$ has $y_a^i$ column vectors which have totally distinct  global phases with  $U_j^i$ for any $0\leq j\leq a-1$. Since for any $a\neq a'\in [t]$, if $p\in R_a$, $q\in R_{a'}$, then $\langle p|U_a^\dagger U_{a'}|q \rangle=0$. Therefore the Eq.(\ref{e2}) can be rewritten as
\begin{equation}
c=\sum\limits_{a=1}^{t-1}\sum\limits_{i=0}^{t-1} y_a^i+x'+(mt+s).\label{e3}
\end{equation}

When $m\geq 3$, because $x'$ runs over $[2,s]$ via Lemma \ref{u2} and  $y_a^i$ runs over $[2,m]$ via Corollary \ref{u4}, so  $c$ can run over  $[mt+s,mt^2+2s]\setminus \{mt+s+1\}$ combining with Lemma \ref{c}.

When $m=2$, we have $y_a^i=2$ via Corollary \ref{u4}, or $y_a^i=0$ by using  $I_2$ instead of $F_{d_i}^{\theta_a}$ for some $1\leq i\leq a$. If $s\geq 3$, since $x'$ runs over $[2,s]$ via Lemma \ref{u2}, then $c$ can run over $[2t+s,2t^2+2s]\setminus \{2t+s+1\}$. If $s=2$, we have $x'=2$ or $x'=0$ by using  $I_2$ instead of $V_s'$, therefore $c$ can run over $[2t+2,2t^2+4]_e$. Consequently, when $m=2$,  there exists a $\rm{QLS}$$(2t+s)$ with cardinality $c$ for any $c\in [2t+s,2t^2+2s]\setminus \{2t+s+1\}$, except possible for $c\in [2t+2,2t^2+4]_o$ with $s=2$. \qed

When $s=1$, we can simplify the Construction \ref{qls1} to obtain the following construction.

\begin{construction}\label{qls2}
If there exists an idempotent $\rm{LS}$$(t)$, an $\rm{LS}$$(m)$ and an $\rm{LS}$$(m+1)$, then there exists a $\rm{QLS}$$(mt+1)$.
\end{construction}

\Proof Suppose $(L(a,b))_{a,b\in[t]}$ is an  idempotent LS($t$). For any $a,b\in [t]$, assume $R_a=[ma, ma+m-1]$ and
\begin{eqnarray*}
\mathcal{A}=\{(a,b,L(a,b))~|~a,b\in [t]\}.
\end{eqnarray*}
 Evidently, $\big(\bigcup_{a\in [t]}R_a\big) \bigcup \{mt\}=[mt+1]$.

 Let $\mathcal{B}=\{|0\rangle, |1\rangle,\ldots,|mt\rangle\}$ be the computational basis of $\mathbb{C}^{mt+1}$.
Assume $W_t=span\{|mt\rangle\}$ and for any $a\in[t]$, let $W_a=span\{|i\rangle: i \in R_a\}.$
Then $\big(\bigoplus\limits_{a\in [t]}W_a\big) \bigoplus W_t=\mathbb{C}^{mt+1}$.
For $a \in [t]$, define unitary matrices of order $mt+1$ as below
\begin{eqnarray}\footnotesize
\centering\small {\begin{tabular}{cc}\scalebox{1}{
$U_a =\begin{pmatrix}
U^0_a& & &\\
&  \ddots & & \\
&  & U^{t-2}_a & \\
& & &V^a_{m+1}
\end{pmatrix}, ~~$ $U_t=
\begin{pmatrix}
I_m& & & \\
&  \ddots & & \\
&  & I_m & \\
& &  & V^t_{m+1}
\end{pmatrix}$},
 \end{tabular}}\label{f4}
 \end{eqnarray}
where $U^i_a$ is a unitary matrix of order $m$, $i\in [t-1]$, $V^t_{m+1}$  a unitary matrix of order $m+1$ and $V^a_{m+1}=\begin{pmatrix}
V^a_m&  \\
&  1
\end{pmatrix}$ with $V^a_m$  a unitary matrix of order $m$.

Next, we construct an array $\Phi=(|\Phi_{x,y}\rangle)$  of order $mt+1$.
For any $A \in \mathcal{A}$, $A=(a,b,L(a,b))$, $a, b\in[t]$, we consider the entries  of $\Phi$ based on  $A$ as follows.

(1) If $a\neq b \in[t]$, since there is an LS($m$), then we construct a QLS($m$), $(|\phi_{i,j}\rangle)_{i,j\in[m]}$  based on the orthonormal basis $\{U_a|mL(a,b)\rangle,U_a|mL(a,b)+1\rangle,\ldots,U_a|mL(a,b)+m-1\rangle\}$ of the subspace $W_{L(a,b)}$. For $i,j\in[m]$, define
\begin{eqnarray}
|\Phi_{ma+i,mb+j}\rangle=|\phi_{i,j}\rangle.
\end{eqnarray}

(2) If $a=b \in [t-1]$, since there is an LS($m+1)$, without loss of generality, let the element $m$ appear  in the lower right corner. Then we can construct a  QLS($m+1$), $(|\varphi_{i,j}\rangle)_{i,j\in[m+1]}$ with $|\varphi_{m,m}\rangle=|mt\rangle$ based on the orthonormal basis $\{U_a|mL(a,b)\rangle,\ldots,U_a|mL(a,b)+m-1\rangle, |mt\rangle\}$ of the subspace $W_{L(a,b)}\bigoplus span\{|mt\rangle\}$. For $i,j\in[m]$, define
\begin{center}
$\begin{array}{l}
|\Phi_{ma+i,mb+j}\rangle=|\varphi_{i,j}\rangle;\\
|\Phi_{ma+i,mt}\rangle=|\varphi_{i,m}\rangle;\\
|\Phi_{mt,mb+j}\rangle=|\varphi_{m,j}\rangle.
\end{array}$
\end{center}

(3) If $a=b=t-1$, since there is an LS($m+1$), then we can construct a  QLS($m+1$), $(|\psi_{i,j}\rangle)_{i,j\in[m]}$ based on the orthonormal basis $\{U_t|mL(t-1,t-1)\rangle, \ldots,U_t|mL(t-1,t-1)+m-1\rangle, U_t|mt\rangle\}$ of the subspace $W_{L(a,b)}\bigoplus span\{|mt\rangle\}$. For any $i,j\in [m+1]$, define
\begin{eqnarray*}
|\Phi_{m(t-1)+i,m(t-1)+j}\rangle=|\psi_{i,j}\rangle.
\end{eqnarray*}

We claim that $(|\Phi_{x,y}\rangle)_{x,y\in[mt+1]}$ is a QLS($mt+1$). In fact, for any $x,y\neq y' \in [mt+1]$, assuming $x=mg+i, y=mh+j,y'=mh'+j'$,  where $0\leq g,h,h'\leq t-1, 0\leq i,j,j'\leq m$, it is necessary to prove the equation below is true,
\begin{eqnarray}
\langle\Phi_{x,y} |\Phi_{x,y'}\rangle=0.\label{55}
\end{eqnarray}

\textcircled{1} For $g\in [t-1]$,  by the cases (1), (2) above,  we need to consider the following cases.

($i$) When $h=h'\in[t], j\neq j'$, then $|\Phi_{x,y}\rangle, |\Phi_{x,y'}\rangle$ are from the same row of $\phi$ or $\varphi$. So the Eq.(\ref{55}) holds.

($ii$) When $h\neq h' \in [t]$, since $L(g,h)\neq L(g,h')$, so $W_{L(g,h)}\cap W_{L(g,h')}=\{0\}$. Moreover $W_{L(g,h)}\cap W_{t}=W_{L(g,h')}\cap W_t=\{0\}$. Consequently, the Eq.(\ref{55}) establishes.

($iii$) When $h=t$ or $h'=t$, without loss of generality, we suppose $h=t$, $y=mt$.  If $h'=g$, because of $\varphi$ being a QLS,  thus the Eq.(\ref{55}) holds.  If $h'\neq g$, so $h\neq h'$. Thus we have $W_{L(g,g)}\cap W_{L(g,h')}=\{0\}$  and  $W_{L(g,h')}\cap W_{t}=\{0\}$,  so  the Eq.(\ref{55}) is true.

\textcircled{2} For $g=t-1$, by the cases (1), (3) above, we need to consider the following cases.

($i$) When $h=h'\in[t]$, then $|\Phi_{x,y}\rangle, |\Phi_{x,y'}\rangle$ are from the same row of $\phi$ or $\psi$. So the Eq.(\ref{55}) holds.

($ii$) When $h\neq h'$, if $h,h'\neq g$, then $L(g,h)\neq L(g,h')$,   $W_{L(g,h)}\cap W_{L(g,h')}=\{0\}$, thus the Eq.(\ref{55}) holds; if $h=g$ or $h'\neq g$, since $W_{L(g,h')}\cap ( W_{L(g,g)} \bigoplus W_{t})= \{0\}$  and $W_{L(g,h)}\cap ( W_{L(g,g)} \bigoplus W_{t})= \{0\}$. Consequently, the Eq.(\ref{55}) establishes.

For any $y, x\neq x'\in[mt+1]$, assuming $y=mh+j$, $x=mg+i$ and $x'=mg'+i'$, where $0\leq g,g',h\leq t-1, 0\leq i,i',j\leq m$, by the same way we can prove
\begin{eqnarray}
\langle \Phi_{x,y}|\Phi_{x',y}\rangle=0.
\end{eqnarray}

Therefore, $\Phi$ is a QLS($mt+1$). \qed

\begin{theorem}\label{c2}
For any $t\geq 3$, $m\geq 2$, there exists a $\rm{QLS}$$(mt+1)$ with cardinality $c$ for any $c\in [mt+1,mt^2+2]\setminus \{mt+2\}$.
\end{theorem}

\Proof We consider the cardinalities of the quantum Latin square $\Phi$ constructed in Construction \ref{qls2}. Suppose $U_a$, $U_t$ are unitary matrices defined in Eq.(\ref{f4}), $a\in[t]$. Specially, we assume $U_0^i=I_m$, $V_{m+1}^0=I_{m+1}$, $i\in[t-1]$. For $1\leq a\leq t-1$, let
\begin{eqnarray}
U_a^i=\begin{pmatrix}
I_{m-d_i}& \\
& F_{d_i}^{\theta_a}
\end{pmatrix}, V_{m}^a=\begin{pmatrix}
I_{m-{d'_a}}& \\
& F_{d'_a}^{\theta'_{a}}
\end{pmatrix}, V_{m+1}^t=\begin{pmatrix}
I_{m+1-k}& \\
& F_k^{\theta}
\end{pmatrix},
\end{eqnarray}
where $F_{d_i}^{\theta_a}$,  $F_{d'_a}^{\theta'_a}$ and $F_k^{\theta}$  are defined in Eq.(\ref{f1}) for $2\leq d_i, d'_a\leq m$, $\theta_a, \theta'_{a}, \theta \in[0,\pi)$, $2\leq k\leq m+1$, $i\in [t-1]$. For $1\leq a\neq a' \leq t-1$, $\theta_a\neq \theta_{a'}$ and $\theta'_a\neq \theta'_{a'}$.

For $1\leq a\leq t-1$, suppose $U_a$ has $x_a$ column vectors which have totally distinct global phases with $U_j$ for any $0\leq j \leq a-1$. Let $V^a_{m+1}$ has $x_a'$ column vectors which have totally distinct global phases with $V^j_{m+1}$ for any $0\leq j \leq a-1$. Assume $z$ is the number of the column vectors in $V_{m+1}^t$ having  totally distinct global phases with $V_{m+1}^a$ for all $a\in[t]$. Thus the cardinality of $\Phi$ is
\begin{equation}
c=\sum\limits_{a=1}^{t-1} x_a+\sum\limits_{a=1}^{t-1}x_a'+z+(mt+1).\label{e3}
\end{equation}
For $i\in[t-1]$, suppose $U_a^i$ has $y_a^i$ column vectors which have totally distinct  global phases with  $U_j^i$ for any $0\leq j\leq a-1$. Since for any $a\neq a'\in [t-1]$, if $p\in R_a$, $q\in R_{a'}$, then $\langle p|U_a^\dagger U_{a'}|q \rangle=0$. Therefore the Eq.(\ref{e3}) can be rewritten as
\begin{equation}
c=\sum\limits_{a=1}^{t-1}\sum\limits_{i=0}^{t-2} y_a^i+\sum\limits_{a=1}^{t-1}x_a'+z+(mt+1).\label{e4}
\end{equation}

Since $m\geq 2$,  then $x'$ runs over $[2,m+1]$,  $y_a^i$ and $x'_a$ runs over $[2,m]$ via Corollary \ref{u4}. Combining with Lemma \ref{c},  the cardinality $c$ of $\Phi$ can run over $[mt+1,mt^2+2]\setminus \{mt+2\}$.  \qed

\subsection{The cardinalities of QLS($mt$)s}
In \cite{Zang2021}, Zang et al. proposed the Direct Product construction for quantum Latin squares, but they did not mention the cardinalities of the obtained QLS($v$)s. To calculate the   cardinality range  of the QLS($v$)s, we simply state the construction process.

\begin{construction}\label{qls3}\rm{\cite{Zang2021}}
$\rm{(Direct~ Product~ construction)}$ If there exists an $\rm{LS}$$(m)$ and an $\rm{LS}$$(t)$, then there exists a $\rm{QLS}$$(mt)$.
\end{construction}
\Proof Let $\mathbb{C}^t=span\{|0\rangle, |1\rangle,\ldots,|t-1\rangle\}$,  $\mathbb{C}^m=span\{|0\rangle, |1\rangle,\ldots,|m-1\rangle\}$ and $\mathbb{C}^{mt}=span\{|0\rangle, |1\rangle,\\ \ldots,|mt-1\rangle\}$, then $\mathbb{C}^{mt}\cong \mathbb{C}^t\otimes \mathbb{C}^m=span\{|i\rangle\otimes |j\rangle~|~i\in[t], j\in[m]\}$.

Since there is an LS($t$) and LS($m$), then we construct a QLS($t$), $(|\phi_{i,j}\rangle)_{i,j\in[t]}$ and a QLS($m$), $(|\varphi_{k,l}\rangle)_{k,l\in[m]}$  based on the computational bases of $\mathbb{C}^t$, $\mathbb{C}^m$ respectively.
For $i \in [t]$, let
\begin{eqnarray}\footnotesize
U_i =\begin{pmatrix}
U^0_i& & &\\
&   U^1_i & & \\
&  & \ddots & \\
& & &U^{t-1}_i
\end{pmatrix}  \label{f5}
\end{eqnarray}
where  $U^j_i$ is a unitary matrix of order $m$ for any  $j\in [t]$. Define $U_0=I_{mt}$.

Suppose $\Phi$ is a square array of order $mt$.
For any $ i,j\in [t]$, $k,l\in[m]$, we define the entries of $\Phi$ as follows
\begin{eqnarray}
|\Phi_{ti+j,mk+l}\rangle=U_i(|\phi_{i,j}\rangle\otimes |\varphi_{k,l}\rangle).
\end{eqnarray}
then $(|\Phi_{ti+j,mk+l}\rangle)_{i,j\in [t], k,l\in[m]}$ is a QLS($mt$). \qed

Specially, we take the submatrix $U_i^j$ of the unitary matrix defined in Eq.(\ref{f5}) as follows
\begin{equation}
U_i^j=\begin{pmatrix}
I_{m-d_j}& \\
& F_{d_j}^{\theta_i}
\end{pmatrix},
\end{equation}
where $F_{d_j}^{\theta_i}$  is the matrix defined in Eq.(\ref{f1}) for $2\leq d_j\leq m$, $\theta_i\neq \theta_{i'}\in[0,\pi)$, for $1\leq i\neq i'\leq t-1$, $j\in [t]$.  Similar to the proof of Theorem \ref{c1} and Theorem \ref{c2}, we have the following results.
\begin{theorem}\label{c3}
For $t\geq 2$, when $m\geq 3$, there exists a $\rm{QLS}$$(mt)$ with cardinality $c$ for any $c\in [mt,mt^2]\setminus \{mt+1\}$; when $m=2$,  there exists a $\rm{QLS}$$(2t)$ with cardinality $c$ for any $c\in [2t,2t^2]_e$.
\end{theorem}

\section{Conclusion}
In this paper, we   use Fourier matrix to construct   QLS($v$)s with maximal cardinality for any $v\geq 4$. For general cardinalities of QLS($v$)s, we respectively used Wilson's construction and Direct Product construction to  establish some possible cardinality range  of a QLS($v$) for any $v\geq 4$. In fact, there are other  good methods for constructing Latin squares in combinatorial designs that can be applied to the study of quantum Latin squares, such as  Filling-holes  construction, PBD construction \cite{Han}. We examine the cardinalities of  QLS($v$)s constructed by these methods and found that although it is still impossible to determine all possible cardinalities,  the range of cardinalities for some special   QLS($v$)s can be significantly expanded. For the sake of the paper's brevity, we  only   construct   a QLS(8) with   cardinality $c\in[8, 64]\backslash\{9\}$, except possible for few values $c\in \{45, 49, 53, 55, 57, 58,59, 61, 62, 63\}$, which is listed in Appendix. It is worth finding new methods to construct QLS($v$)s with all possible cardinalities.

\section*{\textbf{Appendix }}
\footnotesize  {\bf Example:}
There is a QLS$(8)$ with cardinality $c\in[8,64]\backslash \{9\}$, except possible for $c\in \{45, 49, 53, 55, 57, 58,59, 61,$ $ 62, 63\}$.\vspace{0.2cm}

\noindent {\bf Proof}~ By Theorem \ref{c3}, when $m=2, t=4$ we have $c \in [8, 32]_e$; when $m=4, t=2$, we have $c\in [8,16]$. Assume $\Psi_k= ({|\Psi_k\rangle}_{ij})_{i,j\in [8]}$ is a QLS(8) with cardinality $c=k$. When $k\in ([16,32]_o \cup [33,64])\backslash \{45, 49, 53, 55, 57, 58,59, 61, 62, 63\}$, then a QLS(8) with cardinality $k$ can be obtained as below.

\begin{enumerate}[start=1]
\item There exists a QLS(8) with $c$ = 17. \label{17}

\begin{center}
\vspace{0.2cm}\centering\scalebox{0.8}{
$
\Psi_{17} =
\begin{array}{!{\color{black}\vline}c!{\color{black}\vline}c!{\color{black}\vline}c!{\color{black}\vline}c!{\color{red}\vline}c!{\color{black}\vline}c!{\color{black}\vline}c!{\color{black}\vline}c!{\color{black}\vline}}
  \hline
  |0\rangle & |1\rangle & |2\rangle & |3\rangle & |4\rangle & |5\rangle & |6\rangle & |7\rangle\\
  \hline
  |1\rangle & |0\rangle & |3\rangle & |2\rangle & |5\rangle & |4\rangle & |7\rangle & |6\rangle\\
  \hline
  |2\rangle & |3\rangle & |0\rangle & |1\rangle & |6\rangle & |7\rangle & |4\rangle & |5\rangle\\
  \hline
  |3\rangle & |2\rangle & |1\rangle & |0\rangle & |7\rangle & |6\rangle & |5\rangle & |4\rangle\\
  \arrayrulecolor{red}\hline
  |4\rangle & |5\rangle & \color{blue}|6_1\rangle & \color{blue}|7_1\rangle & \color{magenta}|0_4\rangle & \color{magenta}|1_4\rangle & \color{magenta}|2_4\rangle & |3\rangle\\
  \arrayrulecolor{black}\hline
  |5\rangle & |4\rangle & \color{blue}|7_1\rangle & \color{blue}|6_1\rangle & \color{magenta}|1_4\rangle  & \color{magenta}|0_4\rangle & |3\rangle & \color{magenta}|2_4\rangle\\
  \hline
  \color{green}|6_2\rangle & \color{green}|7_2\rangle & \color{cyan}|4_3\rangle & \color{cyan}|5_3\rangle & \color{magenta}|2_4\rangle & |3\rangle & \color{magenta}|0_4\rangle & \color{magenta}|1_4\rangle\\
  \hline
  \color{green}|7_2\rangle & \color{green}|6_2\rangle & \color{cyan}|5_3\rangle & \color{cyan}|4_3\rangle & |3\rangle & \color{magenta}|2_4\rangle & \color{magenta}|1_4\rangle & \color{magenta}|0_4\rangle\\
  \hline
\end{array}
$}
\end{center}

\noindent where
$$ \begin{array}{lll}
\color{blue}|6_1\rangle \color{black}= \frac{1}{7}|6\rangle + \frac{4\sqrt3}{7}|7\rangle,& \color{blue}|7_1\rangle \color{black}=\frac{4\sqrt3}{7}|6\rangle - \frac{1}{7}|7\rangle,&
\color{green}{|6_2\rangle} \color{black}= \frac{1}{2}|6\rangle + \frac{\sqrt3}{2}|7\rangle,\\
\color{green}|7_2\rangle \color{black}=\frac{\sqrt3}{2}|6\rangle - \frac{1}{2}|7\rangle,&
\color{cyan}|4_3\rangle \color{black}= \frac{1}{7}|4\rangle + \frac{4\sqrt3}{7}|5\rangle,&
\color{cyan}|5_3\rangle \color{black}=\frac{4\sqrt3}{7}|4\rangle - \frac{1}{7}|5\rangle,\\
\color{magenta}|0_4\rangle \color{black}= \frac{1}{\sqrt3}(|0\rangle + |1\rangle + |2\rangle),&
\color{magenta}|1_4\rangle \color{black}= \frac{1}{\sqrt3}(|0\rangle + e^{\frac{2\pi \mathrm{i}}{3}}|1\rangle + e^{\frac{\pi \mathrm{i}}{3}}|2\rangle),&
\color{magenta}|2_4\rangle \color{black}= \frac{1}{\sqrt3}(|0\rangle + e^{\frac{\pi \mathrm{i}}{3}}|1\rangle + e^{\frac{2\pi \mathrm{i}}{3}}|2\rangle).
\end{array} $$

\item There exists a QLS(8) with $c$ = 19.\label{19}

\begin{center}
\vspace{0.2cm}\centering\scalebox{0.8}{
$
\Psi_{19} =
\begin{array}{!{\color{black}\vline}c!{\color{black}\vline}c!{\color{black}\vline}c!{\color{black}\vline}c!{\color{red}\vline}c!{\color{black}\vline}c!{\color{black}\vline}c!{\color{black}\vline}c!{\color{black}\vline}}
  \hline
  |0\rangle & |1\rangle & |2\rangle & |3\rangle & |4\rangle & |5\rangle & |6\rangle & |7\rangle\\
  \hline
  |1\rangle & |0\rangle & |3\rangle & |2\rangle & |5\rangle & |4\rangle & |7\rangle & |6\rangle\\
  \hline
  |2\rangle & |3\rangle & |0\rangle & |1\rangle & |6\rangle & |7\rangle & |4\rangle & |5\rangle\\
  \hline
  |3\rangle & |2\rangle & |1\rangle & |0\rangle & |7\rangle & |6\rangle & |5\rangle & |4\rangle\\
  \arrayrulecolor{red}\hline
  \color{blue!50!green}|4_5\rangle & \color{blue!50!green}|5_5\rangle & \color{blue}|6_1\rangle & \color{blue}|7_1\rangle & \color{magenta}|0_4\rangle & \color{magenta}|1_4\rangle & \color{magenta}|2_4\rangle & |3\rangle\\
  \arrayrulecolor{black}\hline
  \color{blue!50!green}|5_5\rangle & \color{blue!50!green}|4_5\rangle & \color{blue}|7_1\rangle & \color{blue}|6_1\rangle & \color{magenta}|1_4\rangle  & \color{magenta}|0_4\rangle & |3\rangle & \color{magenta}|2_4\rangle\\
  \hline
  \color{green}|6_2\rangle & \color{green}|7_2\rangle & \color{cyan}|4_3\rangle & \color{cyan}|5_3\rangle & \color{magenta}|2_4\rangle & |3\rangle & \color{magenta}|0_4\rangle & \color{magenta}|1_4\rangle\\
  \hline
  \color{green}|7_2\rangle & \color{green}|6_2\rangle & \color{cyan}|5_3\rangle & \color{cyan}|4_3\rangle & |3\rangle & \color{magenta}|2_4\rangle & \color{magenta}|1_4\rangle & \color{magenta}|0_4\rangle\\
  \hline
\end{array}
$}
\end{center}

\noindent where
$\color{blue}|6_1\rangle,~\color{blue}|7_1\rangle,~\color{green}|6_2\rangle,~\color{green}|7_2\rangle,
 \color{cyan}|4_3\rangle ,~\color{cyan}|5_3\rangle,~\color{magenta}|0_4\rangle ,~ \color{magenta}|1_4\rangle, \color{magenta}|2_4\rangle$ are given  by Case 1 and
{\footnotesize \begin{eqnarray*}
\color{blue!50!green}|4_5\rangle \color{black} = \frac{1}{2}|4\rangle + \frac{\sqrt3}{2}|5\rangle,~~~~~~\color{blue!50!green}|5_5\rangle \color{black}=\frac{\sqrt3}{2}|4\rangle - \frac{1}{2}|5\rangle.
\end{eqnarray*}}

\item There exists a QLS(8) with $c$ = 21.\label{21}

\begin{center}
\vspace{0.2cm}\centering\scalebox{0.8}{
$
\Psi_{21} =
\begin{array}{!{\color{black}\vline}c!{\color{black}\vline}c!{\color{black}\vline}c!{\color{black}\vline}c!{\color{red}\vline}c!{\color{black}\vline}c!{\color{black}\vline}c!{\color{black}\vline}c!{\color{black}\vline}}
  \hline
  |0\rangle & |1\rangle & |2\rangle & |3\rangle & |4\rangle & |5\rangle & |6\rangle & |7\rangle\\
  \hline
  |1\rangle & |0\rangle & |3\rangle & |2\rangle & |5\rangle & |4\rangle & |7\rangle & |6\rangle\\
  \hline
  |2\rangle & |3\rangle & |0\rangle & |1\rangle & \color{red!50!yellow}|6_6\rangle & \color{red!50!yellow}|7_6\rangle & |4\rangle & |5\rangle\\
  \hline
  |3\rangle & |2\rangle & |1\rangle & |0\rangle & \color{red!50!yellow}|7_6\rangle & \color{red!50!yellow}|6_6\rangle & |5\rangle & |4\rangle\\
  \arrayrulecolor{red}\hline
  \color{blue!50!green}|4_5\rangle & \color{blue!50!green}|5_5\rangle & \color{blue}|6_1\rangle & \color{blue}|7_1\rangle & \color{magenta}|0_4\rangle & \color{magenta}|1_4\rangle & \color{magenta}|2_4\rangle & |3\rangle\\
  \arrayrulecolor{black}\hline
  \color{blue!50!green}|5_5\rangle & \color{blue!50!green}|4_5\rangle & \color{blue}|7_1\rangle & \color{blue}|6_1\rangle & \color{magenta}|1_4\rangle  & \color{magenta}|0_4\rangle & |3\rangle & \color{magenta}|2_4\rangle\\
  \hline
  \color{green}|6_2\rangle & \color{green}|7_2\rangle & \color{cyan}|4_3\rangle & \color{cyan}|5_3\rangle & \color{magenta}|2_4\rangle & |3\rangle & \color{magenta}|0_4\rangle & \color{magenta}|1_4\rangle\\
  \hline
  \color{green}|7_2\rangle & \color{green}|6_2\rangle & \color{cyan}|5_3\rangle & \color{cyan}|4_3\rangle & |3\rangle & \color{magenta}|2_4\rangle & \color{magenta}|1_4\rangle & \color{magenta}|0_4\rangle\\
  \hline
\end{array}
$}
\end{center}

where
$\color{blue}|6_1\rangle,~\color{blue}|7_1\rangle,~\color{green}|6_2\rangle,~\color{green}|7_2\rangle,
 \color{cyan}|4_3\rangle ,~\color{cyan}|5_3\rangle,~\color{magenta}|0_4\rangle ,~ \color{magenta}|1_4\rangle , \color{magenta}|2_4\rangle$ are given  by Case 1,
$\color{blue!50!green}|4_5\rangle ,~ \color{blue!50!green}|5_5\rangle$ are given  by Case 2 and
\begin{eqnarray*}
\color{red!50!yellow}{|6_6\rangle} \color{black}= \frac{1}{8}|6\rangle + \frac{3\sqrt7}{8}|7\rangle,~~~~~~\color{red!50!yellow}|7_6\rangle \color{black}=\frac{3\sqrt7}{8}|6\rangle - \frac{1}{8}|7\rangle.
\end{eqnarray*}

\item There exists a QLS(8) with $c$ = 23.\label{23}

\begin{center}
\vspace{0.2cm}\centering\scalebox{0.8}{
$
\Psi_{23} =
\begin{array}{!{\color{black}\vline}c!{\color{black}\vline}c!{\color{black}\vline}c!{\color{black}\vline}c!{\color{red}\vline}c!{\color{black}\vline}c!{\color{black}\vline}c!{\color{black}\vline}c!{\color{black}\vline}}
  \hline
  |0\rangle & |1\rangle & |2\rangle & |3\rangle & |4\rangle & |5\rangle & |6\rangle & |7\rangle\\
  \hline
  |1\rangle & |0\rangle & |3\rangle & |2\rangle & |5\rangle & |4\rangle & |7\rangle & |6\rangle\\
  \hline
  |2\rangle & |3\rangle & |0\rangle & |1\rangle & \color{red!50!yellow}|6_6\rangle & \color{red!50!yellow}|7_6\rangle & \color{blue!30!green}|4_7\rangle & \color{blue!30!green}|5_7\rangle\\
  \hline
  |3\rangle & |2\rangle & |1\rangle & |0\rangle & \color{red!50!yellow}|7_6\rangle & \color{red!50!yellow}|6_6\rangle & \color{blue!30!green}|5_7\rangle & \color{blue!30!green}|4_7\rangle\\
  \arrayrulecolor{red}\hline
  \color{blue!50!green}|4_5\rangle & \color{blue!50!green}|5_5\rangle & \color{blue}|6_1\rangle & \color{blue}|7_1\rangle & \color{magenta}|0_4\rangle & \color{magenta}|1_4\rangle & \color{magenta}|2_4\rangle & |3\rangle\\
  \arrayrulecolor{black}\hline
  \color{blue!50!green}|5_5\rangle & \color{blue!50!green}|4_5\rangle & \color{blue}|7_1\rangle & \color{blue}|6_1\rangle & \color{magenta}|1_4\rangle  & \color{magenta}|0_4\rangle & |3\rangle & \color{magenta}|2_4\rangle\\
  \hline
  \color{green}|6_2\rangle & \color{green}|7_2\rangle & \color{cyan}|4_3\rangle & \color{cyan}|5_3\rangle & \color{magenta}|2_4\rangle & |3\rangle & \color{magenta}|0_4\rangle & \color{magenta}|1_4\rangle\\
  \hline
  \color{green}|7_2\rangle & \color{green}|6_2\rangle & \color{cyan}|5_3\rangle & \color{cyan}|4_3\rangle & |3\rangle & \color{magenta}|2_4\rangle & \color{magenta}|1_4\rangle & \color{magenta}|0_4\rangle\\
  \hline
\end{array}
$}
\end{center}
where
$\color{blue}|6_1\rangle,~\color{blue}|7_1\rangle,~\color{green}|6_2\rangle,~\color{green}|7_2\rangle,
 \color{cyan}|4_3\rangle ,~\color{cyan}|5_3\rangle,~\color{magenta}|0_4\rangle ,~ \color{magenta}|1_4\rangle , \color{magenta}|2_4\rangle$ are given by Case 1,
$\color{blue!50!green}|4_5\rangle ,~ \color{blue!50!green}|5_5\rangle$, are given by Case 2,
$\color{red!50!yellow}|6_6\rangle ,~  \color{red!50!yellow}|7_6\rangle$ are given by Case 3  and
\begin{eqnarray*}
\color{blue!30!green}{|4_7\rangle} \color{black}= \frac{1}{8}|4\rangle + \frac{3\sqrt7}{8}|5\rangle,~~~~~~\color{blue!30!green}|5_7\rangle \color{black}=\frac{3\sqrt7}{8}|4\rangle - \frac{1}{8}|5\rangle.
\end{eqnarray*}

\item There exists a QLS(8) with $c$ = 25.\label{25}

\begin{center}
\vspace{0.2cm}\centering\scalebox{0.8}{
$
\Psi_{25} =
\begin{array}{!{\color{black}\vline}c!{\color{black}\vline}c!{\color{black}\vline}c!{\color{black}\vline}c!{\color{red}\vline}c!{\color{black}\vline}c!{\color{black}\vline}c!{\color{black}\vline}c!{\color{black}\vline}}
  \hline
  |0\rangle & |1\rangle & |2\rangle & |3\rangle & |4\rangle & |5\rangle & |6\rangle & |7\rangle\\
  \hline
  |1\rangle & |0\rangle & |3\rangle & |2\rangle & |5\rangle & |4\rangle & |7\rangle & |6\rangle\\
  \hline
  \color{red!60!black}|2_8\rangle & \color{red!60!black}|3_8\rangle & |0\rangle & |1\rangle & \color{red!50!yellow}|6_6\rangle & \color{red!50!yellow}|7_6\rangle & \color{blue!30!green}|4_7\rangle & \color{blue!30!green}|5_7\rangle\\
  \hline
  \color{red!60!black}|3_8\rangle & \color{red!60!black}|2_8\rangle & |1\rangle & |0\rangle & \color{red!50!yellow}|7_6\rangle & \color{red!50!yellow}|6_6\rangle & \color{blue!30!green}|5_7\rangle & \color{blue!30!green}|4_7\rangle\\
  \arrayrulecolor{red}\hline
  \color{blue!50!green}|4_5\rangle & \color{blue!50!green}|5_5\rangle & \color{blue}|6_1\rangle & \color{blue}|7_1\rangle & \color{magenta}|0_4\rangle & \color{magenta}|1_4\rangle & \color{magenta}|2_4\rangle & |3\rangle\\
  \arrayrulecolor{black}\hline
  \color{blue!50!green}|5_5\rangle & \color{blue!50!green}|4_5\rangle & \color{blue}|7_1\rangle & \color{blue}|6_1\rangle & \color{magenta}|1_4\rangle  & \color{magenta}|0_4\rangle & |3\rangle & \color{magenta}|2_4\rangle\\
  \hline
  \color{green}|6_2\rangle & \color{green}|7_2\rangle & \color{cyan}|4_3\rangle & \color{cyan}|5_3\rangle & \color{magenta}|2_4\rangle & |3\rangle & \color{magenta}|0_4\rangle & \color{magenta}|1_4\rangle\\
  \hline
  \color{green}|7_2\rangle & \color{green}|6_2\rangle & \color{cyan}|5_3\rangle & \color{cyan}|4_3\rangle & |3\rangle & \color{magenta}|2_4\rangle & \color{magenta}|1_4\rangle & \color{magenta}|0_4\rangle\\
  \hline
\end{array}
$}
\end{center}
where
$\color{blue}|6_1\rangle,~\color{blue}|7_1\rangle,~\color{green}|6_2\rangle,~\color{green}|7_2\rangle,
 \color{cyan}|4_3\rangle ,~\color{cyan}|5_3\rangle,~\color{magenta}|0_4\rangle,~ \color{magenta}|1_4\rangle , \color{magenta}|2_4\rangle$ are given by Case 1,
$\color{blue!50!green}|4_5\rangle ,~ \color{blue!50!green}|5_5\rangle$ are given by Case 2,
$\color{red!50!yellow}|6_6\rangle ,~  \color{red!50!yellow}|7_6\rangle$ are given by Case 3,
$\color{blue!30!green}|4_7\rangle ,~  \color{blue!30!green}|5_7\rangle$ are given by Case 4 and
\begin{eqnarray*}
\color{red!60!black}{|2_8\rangle} \color{black}= \frac{1}{8}|2\rangle + \frac{3\sqrt7}{8}|3\rangle,~~~~~~\color{red!60!black}|3_8\rangle \color{black}=\frac{3\sqrt7}{8}|2\rangle - \frac{1}{8}|3\rangle.
\end{eqnarray*}

\item There exists a QLS(8) with $c$ = 27.\label{27}

\begin{center}
\vspace{0.2cm}\centering\scalebox{0.8}{
$
\Psi_{27} =
\begin{array}{!{\color{black}\vline}c!{\color{black}\vline}c!{\color{black}\vline}c!{\color{black}\vline}c!{\color{red}\vline}c!{\color{black}\vline}c!{\color{black}\vline}c!{\color{black}\vline}c!{\color{black}\vline}}
  \hline
  |0\rangle & |1\rangle & |2\rangle & |3\rangle & |4\rangle & |5\rangle & |6\rangle & |7\rangle\\
  \hline
  |1\rangle & |0\rangle & |3\rangle & |2\rangle & |5\rangle & |4\rangle & |7\rangle & |6\rangle\\
  \hline
  \color{red!60!black}|2_8\rangle & \color{red!60!black}|3_8\rangle & \color{yellow!50!black}|0_9\rangle & \color{yellow!50!black}|1_9\rangle & \color{red!50!yellow}|6_6\rangle & \color{red!50!yellow}|7_6\rangle & \color{blue!30!green}|4_7\rangle & \color{blue!30!green}|5_7\rangle\\
  \hline
  \color{red!60!black}|3_8\rangle & \color{red!60!black}|2_8\rangle & \color{yellow!50!black}|1_9\rangle & \color{yellow!50!black}|0_9\rangle & \color{red!50!yellow}|7_6\rangle & \color{red!50!yellow}|6_6\rangle & \color{blue!30!green}|5_7\rangle & \color{blue!30!green}|4_7\rangle\\
  \arrayrulecolor{red}\hline
  \color{blue!50!green}|4_5\rangle & \color{blue!50!green}|5_5\rangle & \color{blue}|6_1\rangle & \color{blue}|7_1\rangle & \color{magenta}|0_4\rangle & \color{magenta}|1_4\rangle & \color{magenta}|2_4\rangle & |3\rangle\\
  \arrayrulecolor{black}\hline
  \color{blue!50!green}|5_5\rangle & \color{blue!50!green}|4_5\rangle & \color{blue}|7_1\rangle & \color{blue}|6_1\rangle & \color{magenta}|1_4\rangle  & \color{magenta}|0_4\rangle & |3\rangle & \color{magenta}|2_4\rangle\\
  \hline
  \color{green}|6_2\rangle & \color{green}|7_2\rangle & \color{cyan}|4_3\rangle & \color{cyan}|5_3\rangle & \color{magenta}|2_4\rangle & |3\rangle & \color{magenta}|0_4\rangle & \color{magenta}|1_4\rangle\\
  \hline
  \color{green}|7_2\rangle & \color{green}|6_2\rangle & \color{cyan}|5_3\rangle & \color{cyan}|4_3\rangle & |3\rangle & \color{magenta}|2_4\rangle & \color{magenta}|1_4\rangle & \color{magenta}|0_4\rangle\\
  \hline
\end{array}
$}
\end{center}
where
$\color{blue}|6_1\rangle,~\color{blue}|7_1\rangle,~\color{green}|6_2\rangle,~\color{green}|7_2\rangle,
 \color{cyan}|4_3\rangle ,~\color{cyan}|5_3\rangle,~\color{magenta}|0_4\rangle,~ \color{magenta}|1_4\rangle , \color{magenta}|2_4\rangle$ are given by Case 1,
$\color{blue!50!green}|4_5\rangle ,~ \color{blue!50!green}|5_5\rangle$ are given by Case 2,
$\color{red!50!yellow}|6_6\rangle ,~  \color{red!50!yellow}|7_6\rangle$ are given by Case 3,
$\color{blue!30!green}|4_7\rangle ,~  \color{blue!30!green}|5_7\rangle$ are given by Case 4,
$\color{red!60!black}|2_8\rangle ,~ \color{red!60!black}|3_8\rangle $ are given by Case 5 and
\begin{eqnarray*}
\color{yellow!50!black}{|0_9\rangle} \color{black}= \frac{1}{8}|0\rangle + \frac{3\sqrt7}{8}|1\rangle,~~~~~~\color{yellow!50!black}|1_9\rangle \color{black}=\frac{3\sqrt7}{8}|0\rangle - \frac{1}{8}|1\rangle.
\end{eqnarray*}

\item There exists a QLS(8) with $c$ = 29.\label{29}

\begin{center}
\vspace{0.2cm}\centering\scalebox{0.8}{
$
\Psi_{29} =
\begin{array}{!{\color{black}\vline}c!{\color{black}\vline}c!{\color{black}\vline}c!{\color{black}\vline}c!{\color{red}\vline}c!{\color{black}\vline}c!{\color{black}\vline}c!{\color{black}\vline}c!{\color{black}\vline}}
  \hline
  |0\rangle & |1\rangle & |2\rangle & |3\rangle & |4\rangle & |5\rangle & |6\rangle & |7\rangle\\
  \hline
  |1\rangle & |0\rangle & |3\rangle & |2\rangle & |5\rangle & |4\rangle & |7\rangle & |6\rangle\\
  \hline
  \color{red!60!black}|2_8\rangle & \color{red!60!black}|3_8\rangle & |0\rangle & |1\rangle & |6\rangle & |7\rangle & |4\rangle & |5\rangle\\
  \hline
  \color{red!60!black}|3_8\rangle & \color{red!60!black}|2_8\rangle & |1\rangle & |0\rangle & |7\rangle & |6\rangle & |5\rangle & |4\rangle\\
  \arrayrulecolor{red}\hline
  \color{red}|4^0\rangle & \color{red}|5^0\rangle & \color{red}|6^0\rangle & \color{red}|7^0\rangle & \color{magenta}|0_4\rangle & \color{magenta}|1_4\rangle & \color{magenta}|2_4\rangle & |3\rangle\\
  \arrayrulecolor{black}\hline
  \color{red}|5^1\rangle & \color{red}|4^1\rangle & \color{red}|7^1\rangle & \color{red}|6^1\rangle & \color{magenta}|1_4\rangle  & \color{magenta}|0_4\rangle & |3\rangle & \color{magenta}|2_4\rangle\\
  \hline
  \color{red}|6^2\rangle & \color{red}|7^2\rangle & \color{red}|4^2\rangle & \color{red}|5^2\rangle & \color{magenta}|2_4\rangle & |3\rangle & \color{magenta}|0_4\rangle & \color{magenta}|1_4\rangle\\
  \hline
  \color{red}|7^3\rangle & \color{red}|6^3\rangle & \color{red}|5^3\rangle & \color{red}|4^3\rangle & |3\rangle & \color{magenta}|2_4\rangle & \color{magenta}|1_4\rangle & \color{magenta}|0_4\rangle\\
  \hline
\end{array}
$}
\end{center}

where $\color{magenta}|0_4\rangle ,~\color{magenta}|1_4\rangle ,~\color{magenta}|2_4\rangle$ are given by Case 1 and
$\color{red!60!black}|2_8\rangle ,~ \color{red!60!black}|3_8\rangle $ are given by Case 5. Moreover,
let $\Phi = (|\Phi_{ij}\rangle)=(\alpha^i_j)$ is a QLS(4) with maximal cardinality as defined in Theorem \ref{cc}. For $4\leq i\leq 7, j \in [4]$, define ${|\Psi_{29}\rangle}_{ij} = $
$
\begin{pmatrix}
{{\mathbf 0}_{4}}\\
\alpha^{i-4}_j
\end{pmatrix}$, where $\mathbf 0_{4}$ is the zero column vector of order 4, i.e.,

$$\begin{array}{lll}
\color{red}|4^0\rangle \color{black} = \frac{1}{2}(|4\rangle + |5\rangle + |6\rangle + |7\rangle),&
\color{red}|5^0\rangle \color{black} =\frac{1}{2}(|4\rangle - \mathrm{i}|5\rangle - |6\rangle + \mathrm{i}|7\rangle),&
\color{red}|6^0\rangle \color{black} = \frac{1}{2}(|4\rangle - |5\rangle + |6\rangle - |7\rangle),
\\
\color{red}|7^0\rangle \color{black} =\frac{1}{2}(|4\rangle + \mathrm{i}|5\rangle - |6\rangle - \mathrm{i}|7\rangle),&
\color{red}|5^1\rangle \color{black} = \frac{1}{2}(|4\rangle - |5\rangle - \mathrm{i}|6\rangle + \mathrm{i}|7\rangle),&
\color{red}|4^1\rangle \color{black} = \frac{1}{2}(|4\rangle + \mathrm{i}|5\rangle + \mathrm{i}|6\rangle - |7\rangle),\\
\color{red}|7^1\rangle \color{black} = \frac{1}{2}(|4\rangle + |5\rangle - \mathrm{i}|6\rangle - \mathrm{i}|7\rangle),&
\color{red}|6^1\rangle \color{black} = \frac{1}{2}(|4\rangle - \mathrm{i}|5\rangle + \mathrm{i}|6\rangle + |7\rangle),&
\color{red}|6^2\rangle \color{black} = \frac{1}{2}(|4\rangle + |5\rangle - |6\rangle - |7\rangle),\\
\color{red}|7^2\rangle \color{black} = \frac{1}{2}(|4\rangle - \mathrm{i}|5\rangle + |6\rangle - \mathrm{i}|7\rangle)&
\color{red}|4^2\rangle \color{black} = \frac{1}{2}(|4\rangle - |5\rangle - |6\rangle + |7\rangle),&
\color{red}|5^2\rangle \color{black} = \frac{1}{2}(|4\rangle + \mathrm{i}|5\rangle + |6\rangle + \mathrm{i}|7\rangle),\\
\color{red}|7^3\rangle \color{black} = \frac{1}{2}(|4\rangle - |5\rangle + \mathrm{i}|6\rangle - \mathrm{i}|7\rangle),&
\color{red}|6^3\rangle \color{black} = \frac{1}{2}(|4\rangle + \mathrm{i}|5\rangle - \mathrm{i}|6\rangle + |7\rangle),&
\color{red}|5^3\rangle \color{black} = \frac{1}{2}(|4\rangle + |5\rangle + \mathrm{i}|6\rangle + \mathrm{i}|7\rangle),\\
\color{red}|4^3\rangle \color{black} = \frac{1}{2}(|4\rangle - \mathrm{i}|5\rangle - \mathrm{i}|6\rangle - |7\rangle).&&
\end{array}$$

\item There exists a QLS(8) with $c$ = 31.\label{31}

\begin{center}
\vspace{0.2cm}\centering\scalebox{0.8}{
$
\Psi_{31} =
\begin{array}{!{\color{black}\vline}c!{\color{black}\vline}c!{\color{black}\vline}c!{\color{black}\vline}c!{\color{red}\vline}c!{\color{black}\vline}c!{\color{black}\vline}c!{\color{black}\vline}c!{\color{black}\vline}}
  \hline
  |0\rangle & |1\rangle & |2\rangle & |3\rangle & |4\rangle & |5\rangle & |6\rangle & |7\rangle\\
  \hline
  |1\rangle & |0\rangle & |3\rangle & |2\rangle & |5\rangle & |4\rangle & |7\rangle & |6\rangle\\
  \hline
  \color{red!60!black}|2_8\rangle & \color{red!60!black}|3_8\rangle & \color{yellow!50!black}|0_9\rangle & \color{yellow!50!black}|1_9\rangle & |6\rangle & |7\rangle & |4\rangle & |5\rangle\\
  \hline
  \color{red!60!black}|3_8\rangle & \color{red!60!black}|2_8\rangle & \color{yellow!50!black}|1_9\rangle & \color{yellow!50!black}|0_9\rangle & |7\rangle & |6\rangle & |5\rangle & |4\rangle\\
  \arrayrulecolor{red}\hline
  \color{red}|4^0\rangle & \color{red}|5^0\rangle & \color{red}|6^0\rangle & \color{red}|7^0\rangle & \color{magenta}|0_4\rangle & \color{magenta}|1_4\rangle & \color{magenta}|2_4\rangle & |3\rangle\\
  \arrayrulecolor{black}\hline
  \color{red}|5^1\rangle & \color{red}|4^1\rangle & \color{red}|7^1\rangle & \color{red}|6^1\rangle & \color{magenta}|1_4\rangle  & \color{magenta}|0_4\rangle & |3\rangle & \color{magenta}|2_4\rangle\\
  \hline
  \color{red}|6^2\rangle & \color{red}|7^2\rangle & \color{red}|4^2\rangle & \color{red}|5^2\rangle & \color{magenta}|2_4\rangle & |3\rangle & \color{magenta}|0_4\rangle & \color{magenta}|1_4\rangle\\
  \hline
  \color{red}|7^3\rangle & \color{red}|6^3\rangle & \color{red}|5^3\rangle & \color{red}|4^3\rangle & |3\rangle & \color{magenta}|2_4\rangle & \color{magenta}|1_4\rangle & \color{magenta}|0_4\rangle\\
  \hline
\end{array}
$}
\end{center}

where $\color{magenta}|0_4\rangle,~\color{magenta}|1_4\rangle ,~\color{magenta}|2_4\rangle$ are given by Case 1, $\color{red!60!black}|2_8\rangle ,~ \color{red!60!black}|3_8\rangle $ are given by Case 5, $\color{yellow!50!black}|0_9\rangle ,~ \color{yellow!50!black}|1_9\rangle$ are given by Case 6 and  $\color{red}|4^i\rangle,~\color{red}|5^i\rangle,~\color{red}|6^i\rangle,~\color{red}|7^i\rangle$ are given by Case 7 for any $0\leq i \leq 3$,.

\item There exists a QLS(8) with $c$ = 33.\label{33}

\begin{center}
\vspace{0.2cm}\centering\scalebox{0.8}{
$
\Psi_{33} =
\begin{array}{!{\color{black}\vline}c!{\color{black}\vline}c!{\color{black}\vline}c!{\color{black}\vline}c!{\color{red}\vline}c!{\color{black}\vline}c!{\color{black}\vline}c!{\color{black}\vline}c!{\color{black}\vline}}
  \hline
  |0\rangle & |1\rangle & |2\rangle & |3\rangle & |4\rangle & |5\rangle & |6\rangle & |7\rangle\\
  \hline
  |1\rangle & |0\rangle & |3\rangle & |2\rangle & |5\rangle & |4\rangle & |7\rangle & |6\rangle\\
  \hline
  \color{red!60!black}|2_8\rangle & \color{red!60!black}|3_8\rangle & \color{yellow!50!black}|0_9\rangle & \color{yellow!50!black}|1_9\rangle & \color{red!50!yellow}|6_6\rangle & \color{red!50!yellow}|7_6\rangle & |4\rangle & |5\rangle\\
  \hline
  \color{red!60!black}|3_8\rangle & \color{red!60!black}|2_8\rangle & \color{yellow!50!black}|1_9\rangle & \color{yellow!50!black}|0_9\rangle & \color{red!50!yellow}|7_6\rangle & \color{red!50!yellow}|6_6\rangle & |5\rangle & |4\rangle\\
  \arrayrulecolor{red}\hline
  \color{red}|4^0\rangle & \color{red}|5^0\rangle & \color{red}|6^0\rangle & \color{red}|7^0\rangle & \color{magenta}|0_4\rangle & \color{magenta}|1_4\rangle & \color{magenta}|2_4\rangle & |3\rangle\\
  \arrayrulecolor{black}\hline
  \color{red}|5^1\rangle & \color{red}|4^1\rangle & \color{red}|7^1\rangle & \color{red}|6^1\rangle & \color{magenta}|1_4\rangle  & \color{magenta}|0_4\rangle & |3\rangle & \color{magenta}|2_4\rangle\\
  \hline
  \color{red}|6^2\rangle & \color{red}|7^2\rangle & \color{red}|4^2\rangle & \color{red}|5^2\rangle & \color{magenta}|2_4\rangle & |3\rangle & \color{magenta}|0_4\rangle & \color{magenta}|1_4\rangle\\
  \hline
  \color{red}|7^3\rangle & \color{red}|6^3\rangle & \color{red}|5^3\rangle & \color{red}|4^3\rangle & |3\rangle & \color{magenta}|2_4\rangle & \color{magenta}|1_4\rangle & \color{magenta}|0_4\rangle\\
  \hline
\end{array}
$}
\end{center}

where $\color{magenta}|0_4\rangle ,~\color{magenta}|1_4\rangle ,~\color{magenta}|2_4\rangle$ are given by Case 1, $\color{red!50!yellow}|6_6\rangle ,~  \color{red!50!yellow}|7_6\rangle$ are given by Case 3, $\color{red!60!black}|2_8\rangle ,~ \color{red!60!black}|3_8\rangle $ are given by Case 5, $\color{yellow!50!black}|0_9\rangle ,~ \color{yellow!50!black}|1_9\rangle$ are given by Case 6 and  $\color{red}|4^i\rangle,~\color{red}|5^i\rangle,~\color{red}|6^i\rangle,~\color{red}|7^i\rangle$ are given by Case 7 for any $0\leq i \leq 3$.

\item There exists a QLS(8) with $c$ = 34.\label{34}

\begin{center}
\vspace{0.2cm}\centering\scalebox{0.8}{
$
\Psi_{34} =
\begin{array}{!{\color{black}\vline}c!{\color{black}\vline}c!{\color{black}\vline}c!{\color{black}\vline}c!{\color{red}\vline}c!{\color{black}\vline}c!{\color{black}\vline}c!{\color{black}\vline}c!{\color{black}\vline}}
  \hline
  |0\rangle & |1\rangle & |2\rangle & |3\rangle & |4\rangle & |5\rangle & |6\rangle & |7\rangle\\
  \hline
  |1\rangle & |0\rangle & |3\rangle & |2\rangle & |5\rangle & |4\rangle & |7\rangle & |6\rangle\\
  \hline
  \color{red!60!black}|2_8\rangle & \color{red!60!black}|3_8\rangle & \color{yellow!50!black}|0_9\rangle & \color{yellow!50!black}|1_9\rangle & \color{red!50!yellow}|6_6\rangle & \color{red!50!yellow}|7_6\rangle & |4\rangle & |5\rangle\\
  \hline
  \color{red!60!black}|3_8\rangle & \color{red!60!black}|2_8\rangle & \color{yellow!50!black}|1_9\rangle & \color{yellow!50!black}|0_9\rangle & \color{red!50!yellow}|7_6\rangle & \color{red!50!yellow}|6_6\rangle & |5\rangle & |4\rangle\\
  \arrayrulecolor{red}\hline
  \color{red}|4^0\rangle & \color{red}|5^0\rangle & \color{red}|6^0\rangle & \color{red}|7^0\rangle & \color{violet}|0_{10}\rangle & \color{violet}|1_{10}\rangle & \color{violet}|2_{10}\rangle & \color{violet}|3_{10}\rangle\\
  \arrayrulecolor{black}\hline
  \color{red}|5^1\rangle & \color{red}|4^1\rangle & \color{red}|7^1\rangle & \color{red}|6^1\rangle & \color{violet}|1_{10}\rangle  & \color{violet}|0_{10}\rangle & \color{violet}|3_{10}\rangle & \color{violet}|2_{10}\rangle\\
  \hline
  \color{red}|6^2\rangle & \color{red}|7^2\rangle & \color{red}|4^2\rangle & \color{red}|5^2\rangle & \color{violet}|2_{10}\rangle & \color{violet}|3_{10}\rangle & \color{violet}|0_{10}\rangle & \color{violet}|1_{10}\rangle\\
  \hline
  \color{red}|7^3\rangle & \color{red}|6^3\rangle & \color{red}|5^3\rangle & \color{red}|4^3\rangle & \color{violet}|3_{10}\rangle & \color{violet}|2_{10}\rangle & \color{violet}|1_{10}\rangle & \color{violet}|0_{10}\rangle\\
  \hline
\end{array}
$}
\end{center}

where $\color{red!50!yellow}|6_6\rangle ,~  \color{red!50!yellow}|7_6\rangle$ are given by Case 3, $\color{red!60!black}|2_8\rangle ,~ \color{red!60!black}|3_8\rangle $ are given by Case 5, $\color{yellow!50!black}|0_9\rangle ,~ \color{yellow!50!black}|1_9\rangle$ are given by Case 6 and $\color{red}|4^i\rangle,~\color{red}|5^i\rangle,~\color{red}|6^i\rangle,~\color{red}|7^i\rangle$ are given by Case 7 for any $0\leq i \leq 3$.  In addition, let $H=$
$
\begin{pmatrix}
F^{\theta}_4 &  \\
    & I_4
\end{pmatrix}
$,~
$\color{violet}|i_{10}\rangle \color{black}~=~H|i\rangle$,~where $0~\leq i~\leq 3$, $  \theta \in (0, \pi)$.

\item There exists a QLS(8) with $c$ = 35.\label{35}

\begin{center}
\vspace{0.2cm}\centering\scalebox{0.8}{
$
\Psi_{35} =
\begin{array}{!{\color{black}\vline}c!{\color{black}\vline}c!{\color{black}\vline}c!{\color{black}\vline}c!{\color{red}\vline}c!{\color{black}\vline}c!{\color{black}\vline}c!{\color{black}\vline}c!{\color{black}\vline}}
  \hline
  |0\rangle & |1\rangle & |2\rangle & |3\rangle & |4\rangle & |5\rangle & |6\rangle & |7\rangle\\
  \hline
  |1\rangle & |0\rangle & |3\rangle & |2\rangle & |5\rangle & |4\rangle & |7\rangle & |6\rangle\\
  \hline
  \color{red!60!black}|2_8\rangle & \color{red!60!black}|3_8\rangle & \color{yellow!50!black}|0_9\rangle & \color{yellow!50!black}|1_9\rangle & \color{red!50!yellow}|6_6\rangle & \color{red!50!yellow}|7_6\rangle & \color{blue!30!green}|4_7\rangle & \color{blue!30!green}|5_7\rangle\\
  \hline
  \color{red!60!black}|3_8\rangle & \color{red!60!black}|2_8\rangle & \color{yellow!50!black}|1_9\rangle & \color{yellow!50!black}|0_9\rangle & \color{red!50!yellow}|7_6\rangle & \color{red!50!yellow}|6_6\rangle & \color{blue!30!green}|5_7\rangle & \color{blue!30!green}|4_7\rangle\\
  \arrayrulecolor{red}\hline
  \color{red}|4^0\rangle & \color{red}|5^0\rangle & \color{red}|6^0\rangle & \color{red}|7^0\rangle & \color{magenta}|0_4\rangle & \color{magenta}|1_4\rangle & \color{magenta}|2_4\rangle & |3\rangle\\
  \arrayrulecolor{black}\hline
  \color{red}|5^1\rangle & \color{red}|4^1\rangle & \color{red}|7^1\rangle & \color{red}|6^1\rangle & \color{magenta}|1_4\rangle  & \color{magenta}|0_4\rangle & |3\rangle & \color{magenta}|2_4\rangle\\
  \hline
  \color{red}|6^2\rangle & \color{red}|7^2\rangle & \color{red}|4^2\rangle & \color{red}|5^2\rangle & \color{magenta}|2_4\rangle & |3\rangle & \color{magenta}|0_4\rangle & \color{magenta}|1_4\rangle\\
  \hline
  \color{red}|7^3\rangle & \color{red}|6^3\rangle & \color{red}|5^3\rangle & \color{red}|4^3\rangle & |3\rangle & \color{magenta}|2_4\rangle & \color{magenta}|1_4\rangle & \color{magenta}|0_4\rangle\\
  \hline
\end{array}
$}
\end{center}

where
$\color{magenta}|0_4\rangle ,~\color{magenta}|1_4\rangle ,~\color{magenta}|2_4\rangle$ are given by Case 1,
$\color{red!50!yellow}|6_6\rangle ,~  \color{red!50!yellow}|7_6\rangle$ are given by Case 3,
$\color{blue!30!green}|4_7\rangle ,~  \color{blue!30!green}|5_7\rangle$ are given by Case 4,
$\color{red!60!black}|2_8\rangle ,~ \color{red!60!black}|3_8\rangle $ are given by Case 5,
$\color{yellow!50!black}|0_9\rangle ,~ \color{yellow!50!black}|1_9\rangle$ are given by Case 6 and
 $\color{red}|4^i\rangle,~\color{red}|5^i\rangle,~\color{red}|6^i\rangle,~\color{red}|7^i\rangle$ are given by Case 7 for any $0\leq i \leq 3$.

\item There exists a QLS(8) with $c$ = 36.\label{36}

\begin{center}
\vspace{0.2cm}\centering\scalebox{0.8}{
$
\Psi_{36} =
\begin{array}{!{\color{black}\vline}c!{\color{black}\vline}c!{\color{black}\vline}c!{\color{black}\vline}c!{\color{red}\vline}c!{\color{black}\vline}c!{\color{black}\vline}c!{\color{black}\vline}c!{\color{black}\vline}}
  \hline
  |0\rangle & |1\rangle & |2\rangle & |3\rangle & |4\rangle & |5\rangle & |6\rangle & |7\rangle\\
  \hline
  |1\rangle & |0\rangle & |3\rangle & |2\rangle & |5\rangle & |4\rangle & |7\rangle & |6\rangle\\
  \hline
  \color{red!60!black}|2_8\rangle & \color{red!60!black}|3_8\rangle & \color{yellow!50!black}|0_9\rangle & \color{yellow!50!black}|1_9\rangle & \color{red!50!yellow}|6_6\rangle & \color{red!50!yellow}|7_6\rangle & \color{blue!30!green}|4_7\rangle & \color{blue!30!green}|5_7\rangle\\
  \hline
  \color{red!60!black}|3_8\rangle & \color{red!60!black}|2_8\rangle & \color{yellow!50!black}|1_9\rangle & \color{yellow!50!black}|0_9\rangle & \color{red!50!yellow}|7_6\rangle & \color{red!50!yellow}|6_6\rangle & \color{blue!30!green}|5_7\rangle & \color{blue!30!green}|4_7\rangle\\
  \arrayrulecolor{red}\hline
  \color{red}|4^0\rangle & \color{red}|5^0\rangle & \color{red}|6^0\rangle & \color{red}|7^0\rangle & \color{violet}|0_{10}\rangle & \color{violet}|1_{10}\rangle & \color{violet}|2_{10}\rangle & \color{violet}|3_{10}\rangle\\
  \arrayrulecolor{black}\hline
  \color{red}|5^1\rangle & \color{red}|4^1\rangle & \color{red}|7^1\rangle & \color{red}|6^1\rangle & \color{violet}|1_{10}\rangle  & \color{violet}|0_{10}\rangle & \color{violet}|3_{10}\rangle & \color{violet}|2_{10}\rangle\\
  \hline
  \color{red}|6^2\rangle & \color{red}|7^2\rangle & \color{red}|4^2\rangle & \color{red}|5^2\rangle & \color{violet}|2_{10}\rangle & \color{violet}|3_{10}\rangle & \color{violet}|0_{10}\rangle & \color{violet}|1_{10}\rangle\\
  \hline
  \color{red}|7^3\rangle & \color{red}|6^3\rangle & \color{red}|5^3\rangle & \color{red}|4^3\rangle & \color{violet}|3_{10}\rangle & \color{violet}|2_{10}\rangle & \color{violet}|1_{10}\rangle & \color{violet}|0_{10}\rangle\\
  \hline
\end{array}
$}
\end{center}

where
$\color{red!50!yellow}|6_6\rangle ,~  \color{red!50!yellow}|7_6\rangle$ are given by Case 3,
$\color{blue!30!green}|4_7\rangle ,~  \color{blue!30!green}|5_7\rangle$ are given by Case 4,
$\color{red!60!black}|2_8\rangle ,~ \color{red!60!black}|3_8\rangle $ are given by Case 5,
$\color{yellow!50!black}|0_9\rangle ,~ \color{yellow!50!black}|1_9\rangle$ are given by Case 6,
$\color{violet}|0_{10}\rangle ,~ \color{violet}|1_{10}\rangle ,~ \color{violet}|2_{10}\rangle ,~ \color{violet}|3_{10}\rangle$ are given by Case 10 and
$\color{red}|4^i\rangle,~\color{red}|5^i\rangle,~\color{red}|6^i\rangle,~\color{red}|7^i\rangle$ are given by Case 7  for any $0\leq i \leq 3$.

\item There exists a QLS(8) with $c$ = 37.\label{37}

\begin{center}
\vspace{0.2cm}\centering\scalebox{0.8}{
$
\Psi_{37} =
\begin{array}{!{\color{black}\vline}c!{\color{black}\vline}c!{\color{black}\vline}c!{\color{black}\vline}c!{\color{red}\vline}c!{\color{black}\vline}c!{\color{black}\vline}c!{\color{black}\vline}c!{\color{black}\vline}}
  \hline
  |0\rangle & |1\rangle & \frac{1}{\sqrt2}(|2\rangle-|3\rangle) & \frac{1}{\sqrt2}(|2\rangle+|3\rangle) & |4\rangle & |5\rangle & |6\rangle & |7\rangle\\
  \hline
  |2\rangle & |3\rangle & \frac{1}{\sqrt2}(|0\rangle-|1\rangle) & \frac{1}{\sqrt2}(|0\rangle+|1\rangle) & |5\rangle & |4\rangle & |7\rangle & |6\rangle\\
  \hline
  \frac{1}{\sqrt2}(|1\rangle+|3\rangle) & \frac{1}{\sqrt2}(|0\rangle-|2\rangle) & \frac{1}{2}(|0\rangle+|1\rangle+|2\rangle-|3\rangle) & \frac{1}{2}(|0\rangle-|1\rangle+|2\rangle+|3\rangle) & \color{red!50!yellow}|6_6\rangle & \color{red!50!yellow}|7_6\rangle & |4\rangle & |5\rangle\\
  \hline
  \frac{1}{\sqrt2}(|1\rangle-|3\rangle) & \frac{1}{\sqrt2}(|0\rangle+|2\rangle) & \frac{1}{2}(|0\rangle+|1\rangle-|2\rangle+|3\rangle) & \frac{1}{2}(|0\rangle-|1\rangle-|2\rangle-|3\rangle) & \color{red!50!yellow}|7_6\rangle & \color{red!50!yellow}|6_6\rangle & |5\rangle & |4\rangle\\
  \arrayrulecolor{red}\hline
  |4\rangle & |5\rangle & \frac{1}{\sqrt2}(|6\rangle-|7\rangle) & \frac{1}{\sqrt2}(|6\rangle+|7\rangle) & \color{magenta}|0_4\rangle & \color{magenta}|1_4\rangle & \color{magenta}|2_4\rangle & |3\rangle\\
  \arrayrulecolor{black}\hline
  |6\rangle & |7\rangle & \frac{1}{\sqrt2}(|0\rangle-|1\rangle) & \frac{1}{\sqrt2}(|4\rangle+|5\rangle) & \color{magenta}|1_4\rangle  & \color{magenta}|0_4\rangle & |3\rangle & \color{magenta}|2_4\rangle\\
  \hline
   \frac{1}{\sqrt2}(|5\rangle+|7\rangle) & \frac{1}{\sqrt2}(|4\rangle-|6\rangle) & \frac{1}{2}(|4\rangle+|5\rangle+|6\rangle-|7\rangle) & \frac{1}{2}(|4\rangle-|5\rangle+|6\rangle+|7\rangle) & \color{magenta}|2_4\rangle & |3\rangle & \color{magenta}|0_4\rangle & \color{magenta}|1_4\rangle\\
  \hline
  \frac{1}{\sqrt2}(|5\rangle-|7\rangle) & \frac{1}{\sqrt2}(|4\rangle+|6\rangle) & \frac{1}{2}(|4\rangle+|5\rangle-|6\rangle+|7\rangle) & \frac{1}{2}(|4\rangle-|5\rangle-|6\rangle-|7\rangle) & |3\rangle & \color{magenta}|2_4\rangle & \color{magenta}|1_4\rangle & \color{magenta}|0_4\rangle\\
  \hline
\end{array}
$}
\end{center}

where
$\color{magenta}|0_4\rangle ,~\color{magenta}|1_4\rangle ,~\color{magenta}|2_4\rangle$ are given by Case 1 and
$\color{red!50!yellow}|6_6\rangle ,~  \color{red!50!yellow}|7_6\rangle$ are given by Case 3.

\item There exists a QLS(8) with $c$ = 38.\label{38}

\begin{center}
\vspace{0.2cm}\centering\scalebox{0.8}{
$
\Psi_{38} =
\begin{array}{!{\color{black}\vline}c!{\color{black}\vline}c!{\color{black}\vline}c!{\color{black}\vline}c!{\color{red}\vline}c!{\color{black}\vline}c!{\color{black}\vline}c!{\color{black}\vline}c!{\color{black}\vline}}
  \hline
  |0\rangle & |1\rangle & \frac{1}{\sqrt2}(|2\rangle-|3\rangle) & \frac{1}{\sqrt2}(|2\rangle+|3\rangle) & |4\rangle & |5\rangle & |6\rangle & |7\rangle\\
  \hline
  |2\rangle & |3\rangle & \frac{1}{\sqrt2}(|0\rangle-|1\rangle) & \frac{1}{\sqrt2}(|0\rangle+|1\rangle) & |5\rangle & |4\rangle & |7\rangle & |6\rangle\\
  \hline
  \frac{1}{\sqrt2}(|1\rangle+|3\rangle) & \frac{1}{\sqrt2}(|0\rangle-|2\rangle) & \frac{1}{2}(|0\rangle+|1\rangle+|2\rangle-|3\rangle) & \frac{1}{2}(|0\rangle-|1\rangle+|2\rangle+|3\rangle) & \color{red!50!yellow}|6_6\rangle & \color{red!50!yellow}|7_6\rangle & |4\rangle & |5\rangle\\
  \hline
  \frac{1}{\sqrt2}(|1\rangle-|3\rangle) & \frac{1}{\sqrt2}(|0\rangle+|2\rangle) & \frac{1}{2}(|0\rangle+|1\rangle-|2\rangle+|3\rangle) & \frac{1}{2}(|0\rangle-|1\rangle-|2\rangle-|3\rangle) & \color{red!50!yellow}|7_6\rangle & \color{red!50!yellow}|6_6\rangle & |5\rangle & |4\rangle\\
  \arrayrulecolor{red}\hline
  |4\rangle & |5\rangle & \frac{1}{\sqrt2}(|6\rangle-|7\rangle) & \frac{1}{\sqrt2}(|6\rangle+|7\rangle) & \color{violet}|0_{10}\rangle & \color{violet}|1_{10}\rangle & \color{violet}|2_{10}\rangle & \color{violet}|3_{10}\rangle\\
  \arrayrulecolor{black}\hline
  |6\rangle & |7\rangle & \frac{1}{\sqrt2}(|0\rangle-|1\rangle) & \frac{1}{\sqrt2}(|4\rangle+|5\rangle) & \color{violet}|1_{10}\rangle  & \color{violet}|0_{10}\rangle & \color{violet}|3_{10}\rangle & \color{violet}|2_{10}\rangle\\
  \hline
  \frac{1}{\sqrt2}(|5\rangle+|7\rangle) & \frac{1}{\sqrt2}(|4\rangle-|6\rangle) & \frac{1}{2}(|4\rangle+|5\rangle+|6\rangle-|7\rangle) & \frac{1}{2}(|4\rangle-|5\rangle+|6\rangle+|7\rangle) & \color{violet}|2_{10}\rangle & \color{violet}|3_{10}\rangle & \color{violet}|0_{10}\rangle & \color{violet}|1_{10}\rangle\\
  \hline
  \frac{1}{\sqrt2}(|5\rangle-|7\rangle) & \frac{1}{\sqrt2}(|4\rangle+|6\rangle) & \frac{1}{2}(|4\rangle+|5\rangle-|6\rangle+|7\rangle) & \frac{1}{2}(|4\rangle-|5\rangle-|6\rangle-|7\rangle) & \color{violet}|3_{10}\rangle & \color{violet}|2_{10}\rangle & \color{violet}|1_{10}\rangle & \color{violet}|0_{10}\rangle\\
  \hline
\end{array}
$}
\end{center}

where $\color{red!50!yellow}|6_6\rangle ,~  \color{red!50!yellow}|7_6\rangle$ are given by Case 3 and
$\color{violet}|0_{10}\rangle ,~ \color{violet}|1_{10}\rangle ,~ \color{violet}|2_{10}\rangle ,~ \color{violet}|3_{10}\rangle$ are given by Case 10.

\item There exists a QLS(8) with $c$ = 39.\label{39}

\begin{center}
\vspace{0.2cm}\centering\scalebox{0.8}{
$
\Psi_{39} =
\begin{array}{!{\color{black}\vline}c!{\color{black}\vline}c!{\color{black}\vline}c!{\color{black}\vline}c!{\color{red}\vline}c!{\color{black}\vline}c!{\color{black}\vline}c!{\color{black}\vline}c!{\color{black}\vline}}
  \hline
  |0\rangle & |1\rangle & \frac{1}{\sqrt2}(|2\rangle-|3\rangle) & \frac{1}{\sqrt2}(|2\rangle+|3\rangle) & |4\rangle & |5\rangle & |6\rangle & |7\rangle\\
  \hline
  |2\rangle & |3\rangle & \frac{1}{\sqrt2}(|0\rangle-|1\rangle) & \frac{1}{\sqrt2}(|0\rangle+|1\rangle) & |5\rangle & |4\rangle & |7\rangle & |6\rangle\\
  \hline
  \frac{1}{\sqrt2}(|1\rangle+|3\rangle) & \frac{1}{\sqrt2}(|0\rangle-|2\rangle) & \frac{1}{2}(|0\rangle+|1\rangle+|2\rangle-|3\rangle) & \frac{1}{2}(|0\rangle-|1\rangle+|2\rangle+|3\rangle) & \color{red!50!yellow}|6_6\rangle & \color{red!50!yellow}|7_6\rangle & \color{blue!30!green}|4_7\rangle & \color{blue!30!green}|5_7\rangle\\
  \hline
  \frac{1}{\sqrt2}(|1\rangle-|3\rangle) & \frac{1}{\sqrt2}(|0\rangle+|2\rangle) & \frac{1}{2}(|0\rangle+|1\rangle-|2\rangle+|3\rangle) & \frac{1}{2}(|0\rangle-|1\rangle-|2\rangle-|3\rangle) & \color{red!50!yellow}|7_6\rangle & \color{red!50!yellow}|6_6\rangle & \color{blue!30!green}|5_7\rangle & \color{blue!30!green}|4_7\rangle\\
  \arrayrulecolor{red}\hline
  |4\rangle & |5\rangle & \frac{1}{\sqrt2}(|6\rangle-|7\rangle) & \frac{1}{\sqrt2}(|6\rangle+|7\rangle) & \color{magenta}|0_4\rangle & \color{magenta}|1_4\rangle & \color{magenta}|2_4\rangle & |3\rangle\\
  \arrayrulecolor{black}\hline
  |6\rangle & |7\rangle & \frac{1}{\sqrt2}(|0\rangle-|1\rangle) & \frac{1}{\sqrt2}(|4\rangle+|5\rangle) & \color{magenta}|1_4\rangle  & \color{magenta}|0_4\rangle & |3\rangle & \color{magenta}|2_4\rangle\\
  \hline
  \frac{1}{\sqrt2}(|5\rangle+|7\rangle) & \frac{1}{\sqrt2}(|4\rangle-|6\rangle) & \frac{1}{2}(|4\rangle+|5\rangle+|6\rangle-|7\rangle) & \frac{1}{2}(|4\rangle-|5\rangle+|6\rangle+|7\rangle) & \color{magenta}|2_4\rangle & |3\rangle & \color{magenta}|0_4\rangle & \color{magenta}|1_4\rangle\\
  \hline
  \frac{1}{\sqrt2}(|5\rangle-|7\rangle) & \frac{1}{\sqrt2}(|4\rangle+|6\rangle) & \frac{1}{2}(|4\rangle+|5\rangle-|6\rangle+|7\rangle) & \frac{1}{2}(|4\rangle-|5\rangle-|6\rangle-|7\rangle) & |3\rangle & \color{magenta}|2_4\rangle & \color{magenta}|1_4\rangle & \color{magenta}|0_4\rangle\\
  \hline
\end{array}
$}
\end{center}

where
$\color{magenta}|0_4\rangle ,~\color{magenta}|1_4\rangle ,~\color{magenta}|2_4\rangle$ are given by Case 1,
$\color{red!50!yellow}|6_6\rangle ,~\color{red!50!yellow}|7_6\rangle$ are given by Case 3 and
$\color{blue!30!green}|4_7\rangle ,~\color{blue!30!green}|5_7\rangle$ are given by Case 4.

\item There exists a QLS(8) with $c$ = 40.\label{40}

\begin{center}
\vspace{0.2cm}\centering\scalebox{0.8}{
$
\Psi_{40} =
\begin{array}{!{\color{black}\vline}c!{\color{black}\vline}c!{\color{black}\vline}c!{\color{black}\vline}c!{\color{red}\vline}c!{\color{black}\vline}c!{\color{black}\vline}c!{\color{black}\vline}c!{\color{black}\vline}}
  \hline
  |0\rangle & |1\rangle & \frac{1}{\sqrt2}(|2\rangle-|3\rangle) & \frac{1}{\sqrt2}(|2\rangle+|3\rangle) & |4\rangle & |5\rangle & |6\rangle & |7\rangle\\
  \hline
  |2\rangle & |3\rangle & \frac{1}{\sqrt2}(|0\rangle-|1\rangle) & \frac{1}{\sqrt2}(|0\rangle+|1\rangle) & |5\rangle & |4\rangle & |7\rangle & |6\rangle\\
  \hline
  \frac{1}{\sqrt2}(|1\rangle+|3\rangle) & \frac{1}{\sqrt2}(|0\rangle-|2\rangle) & \frac{1}{2}(|0\rangle+|1\rangle+|2\rangle-|3\rangle) & \frac{1}{2}(|0\rangle-|1\rangle+|2\rangle+|3\rangle) & \color{red!50!yellow}|6_6\rangle & \color{red!50!yellow}|7_6\rangle & \color{blue!30!green}|4_7\rangle & \color{blue!30!green}|5_7\rangle\\
  \hline
  \frac{1}{\sqrt2}(|1\rangle-|3\rangle) & \frac{1}{\sqrt2}(|0\rangle+|2\rangle) & \frac{1}{2}(|0\rangle+|1\rangle-|2\rangle+|3\rangle) & \frac{1}{2}(|0\rangle-|1\rangle-|2\rangle-|3\rangle) & \color{red!50!yellow}|7_6\rangle & \color{red!50!yellow}|6_6\rangle & \color{blue!30!green}|5_7\rangle & \color{blue!30!green}|4_7\rangle\\
  \arrayrulecolor{red}\hline
  |4\rangle & |5\rangle & \frac{1}{\sqrt2}(|6\rangle-|7\rangle) & \frac{1}{\sqrt2}(|6\rangle+|7\rangle) & \color{violet}|0_{10}\rangle & \color{violet}|1_{10}\rangle & \color{violet}|2_{10}\rangle & \color{violet}|3_{10}\rangle\\
  \arrayrulecolor{black}\hline
  |6\rangle & |7\rangle & \frac{1}{\sqrt2}(|0\rangle-|1\rangle) & \frac{1}{\sqrt2}(|4\rangle+|5\rangle) & \color{violet}|1_{10}\rangle  & \color{violet}|0_{10}\rangle & \color{violet}|3_{10}\rangle & \color{violet}|2_{10}\rangle\\
  \hline
  \frac{1}{\sqrt2}(|5\rangle+|7\rangle) & \frac{1}{\sqrt2}(|4\rangle-|6\rangle) & \frac{1}{2}(|4\rangle+|5\rangle+|6\rangle-|7\rangle) & \frac{1}{2}(|4\rangle-|5\rangle+|6\rangle+|7\rangle) & \color{violet}|2_{10}\rangle & \color{violet}|3_{10}\rangle & \color{violet}|0_{10}\rangle & \color{violet}|1_{10}\rangle\\
  \hline
  \frac{1}{\sqrt2}(|5\rangle-|7\rangle) & \frac{1}{\sqrt2}(|4\rangle+|6\rangle) & \frac{1}{2}(|4\rangle+|5\rangle-|6\rangle+|7\rangle) & \frac{1}{2}(|4\rangle-|5\rangle-|6\rangle-|7\rangle) & \color{violet}|3_{10}\rangle & \color{violet}|2_{10}\rangle & \color{violet}|1_{10}\rangle & \color{violet}|0_{10}\rangle\\
  \hline
\end{array}
$}
\end{center}

where $\color{red!50!yellow}|6_6\rangle ,~  \color{red!50!yellow}|7_6\rangle$ are given by Case 3,
$\color{blue!30!green}|4_7\rangle ,~\color{blue!30!green}|5_7\rangle$ are given by Case 4 and
$\color{violet}|0_{10}\rangle ,~ \color{violet}|1_{10}\rangle ,~ \color{violet}|2_{10}\rangle ,~ \color{violet}|3_{10}\rangle$ are given by Case 10.

\item There exists a QLS(8) with $c$ = 41.\label{41}

\begin{center}
\vspace{0.2cm}\centering\scalebox{0.8}{
$
\Psi_{41} =
\begin{array}{!{\color{black}\vline}c!{\color{black}\vline}c!{\color{black}\vline}c!{\color{black}\vline}c!{\color{red}\vline}c!{\color{black}\vline}c!{\color{black}\vline}c!{\color{black}\vline}c!{\color{black}\vline}}
  \hline
  |0\rangle & |1\rangle & \frac{1}{\sqrt2}(|2\rangle-|3\rangle) & \frac{1}{\sqrt2}(|2\rangle+|3\rangle) & |4\rangle & |5\rangle & |6\rangle & |7\rangle\\
  \hline
  |2\rangle & |3\rangle & \frac{1}{\sqrt2}(|0\rangle-|1\rangle) & \frac{1}{\sqrt2}(|0\rangle+|1\rangle) & |5\rangle & |4\rangle & |7\rangle & |6\rangle\\
  \hline
  \frac{1}{\sqrt2}(|1\rangle+|3\rangle) & \frac{1}{\sqrt2}(|0\rangle-|2\rangle) & \frac{1}{2}(|0\rangle+|1\rangle+|2\rangle-|3\rangle) & \frac{1}{2}(|0\rangle-|1\rangle+|2\rangle+|3\rangle) & \color{red!50!yellow}|6_6\rangle & \color{red!50!yellow}|7_6\rangle & |4\rangle & |5\rangle\\
  \hline
  \frac{1}{\sqrt2}(|1\rangle-|3\rangle) & \frac{1}{\sqrt2}(|0\rangle+|2\rangle) & \frac{1}{2}(|0\rangle+|1\rangle-|2\rangle+|3\rangle) & \frac{1}{2}(|0\rangle-|1\rangle-|2\rangle-|3\rangle) & \color{red!50!yellow}|7_6\rangle & \color{red!50!yellow}|6_6\rangle & |5\rangle & |4\rangle\\
  \arrayrulecolor{red}\hline
  \color{red}|4^0\rangle & \color{red}|5^0\rangle & \color{red}|6^0\rangle & \color{red}|7^0\rangle & \color{magenta}|0_4\rangle & \color{magenta}|1_4\rangle & \color{magenta}|2_4\rangle & |3\rangle\\
  \arrayrulecolor{black}\hline
  \color{red}|5^1\rangle & \color{red}|4^1\rangle & \color{red}|7^1\rangle & \color{red}|6^1\rangle & \color{magenta}|1_4\rangle  & \color{magenta}|0_4\rangle & |3\rangle & \color{magenta}|2_4\rangle\\
  \hline
  \color{red}|6^2\rangle & \color{red}|7^2\rangle & \color{red}|4^2\rangle & \color{red}|5^2\rangle & \color{magenta}|2_4\rangle & |3\rangle & \color{magenta}|0_4\rangle & \color{magenta}|1_4\rangle\\
  \hline
  \color{red}|7^3\rangle & \color{red}|6^3\rangle & \color{red}|5^3\rangle & \color{red}|4^3\rangle & |3\rangle & \color{magenta}|2_4\rangle & \color{magenta}|1_4\rangle & \color{magenta}|0_4\rangle\\
  \hline
\end{array}
$}
\end{center}

where
$\color{magenta}|0_4\rangle ,~\color{magenta}|1_4\rangle ,~\color{magenta}|2_4\rangle$ are given by Case 1,
$\color{red!50!yellow}|6_6\rangle ,~  \color{red!50!yellow}|7_6\rangle$ are given by Case 3 and
 $\color{red}|4^i\rangle,~\color{red}|5^i\rangle,~\color{red}|6^i\rangle,~\color{red}|7^i\rangle$ are given by Case 7 for  $0\leq i \leq 3$.

\item There exists a QLS(8) with $c$ = 42.\label{42}

\begin{center}
\vspace{0.2cm}\centering\scalebox{0.8}{
$
\Psi_{42} =
\begin{array}{!{\color{black}\vline}c!{\color{black}\vline}c!{\color{black}\vline}c!{\color{black}\vline}c!{\color{red}\vline}c!{\color{black}\vline}c!{\color{black}\vline}c!{\color{black}\vline}c!{\color{black}\vline}}
  \hline
  |0\rangle & |1\rangle & \frac{1}{\sqrt2}(|2\rangle-|3\rangle) & \frac{1}{\sqrt2}(|2\rangle+|3\rangle) & |4\rangle & |5\rangle & |6\rangle & |7\rangle\\
  \hline
  |2\rangle & |3\rangle & \frac{1}{\sqrt2}(|0\rangle-|1\rangle) & \frac{1}{\sqrt2}(|0\rangle+|1\rangle) & |5\rangle & |4\rangle & |7\rangle & |6\rangle\\
  \hline
  \frac{1}{\sqrt2}(|1\rangle+|3\rangle) & \frac{1}{\sqrt2}(|0\rangle-|2\rangle) & \frac{1}{2}(|0\rangle+|1\rangle+|2\rangle-|3\rangle) & \frac{1}{2}(|0\rangle-|1\rangle+|2\rangle+|3\rangle) & \color{red!50!yellow}|6_6\rangle & \color{red!50!yellow}|7_6\rangle & |4\rangle & |5\rangle\\
  \hline
  \frac{1}{\sqrt2}(|1\rangle-|3\rangle) & \frac{1}{\sqrt2}(|0\rangle+|2\rangle) & \frac{1}{2}(|0\rangle+|1\rangle-|2\rangle+|3\rangle) & \frac{1}{2}(|0\rangle-|1\rangle-|2\rangle-|3\rangle) & \color{red!50!yellow}|7_6\rangle & \color{red!50!yellow}|6_6\rangle & |5\rangle & |4\rangle\\
  \arrayrulecolor{red}\hline
  \color{red}|4^0\rangle & \color{red}|5^0\rangle & \color{red}|6^0\rangle & \color{red}|7^0\rangle & \color{violet}|0_{10}\rangle & \color{violet}|1_{10}\rangle & \color{violet}|2_{10}\rangle & \color{violet}|3_{10}\rangle\\
  \arrayrulecolor{black}\hline
  \color{red}|5^1\rangle & \color{red}|4^1\rangle & \color{red}|7^1\rangle & \color{red}|6^1\rangle & \color{violet}|1_{10}\rangle  & \color{violet}|0_{10}\rangle & \color{violet}|3_{10}\rangle & \color{violet}|2_{10}\rangle\\
  \hline
  \color{red}|6^2\rangle & \color{red}|7^2\rangle & \color{red}|4^2\rangle & \color{red}|5^2\rangle & \color{violet}|2_{10}\rangle & \color{violet}|3_{10}\rangle & \color{violet}|0_{10}\rangle & \color{violet}|1_{10}\rangle\\
  \hline
  \color{red}|7^3\rangle & \color{red}|6^3\rangle & \color{red}|5^3\rangle & \color{red}|4^3\rangle & \color{violet}|3_{10}\rangle & \color{violet}|2_{10}\rangle & \color{violet}|1_{10}\rangle & \color{violet}|0_{10}\rangle\\
  \hline
\end{array}
$}
\end{center}

where
$\color{red!50!yellow}|6_6\rangle ,~  \color{red!50!yellow}|7_6\rangle$ are given by Case 3,
$\color{violet}|0_{10}\rangle ,~ \color{violet}|1_{10}\rangle ,~ \color{violet}|2_{10}\rangle ,~ \color{violet}|3_{10}\rangle$ are given by Case 10 and
$\color{red}|4^i\rangle,~\color{red}|5^i\rangle,~\color{red}|6^i\rangle,~\color{red}|7^i\rangle$ are given by Case 7 for any $0\leq i \leq 3$, .

\item There exists a QLS(8) with $c$ = 43.\label{43}

\begin{center}
\vspace{0.2cm}\centering\scalebox{0.8}{
$
\Psi_{43} =
\begin{array}{!{\color{black}\vline}c!{\color{black}\vline}c!{\color{black}\vline}c!{\color{black}\vline}c!{\color{red}\vline}c!{\color{black}\vline}c!{\color{black}\vline}c!{\color{black}\vline}c!{\color{black}\vline}}
  \hline
  |0\rangle & |1\rangle & \frac{1}{\sqrt2}(|2\rangle-|3\rangle) & \frac{1}{\sqrt2}(|2\rangle+|3\rangle) & |4\rangle & |5\rangle & |6\rangle & |7\rangle\\
  \hline
  |2\rangle & |3\rangle & \frac{1}{\sqrt2}(|0\rangle-|1\rangle) & \frac{1}{\sqrt2}(|0\rangle+|1\rangle) & |5\rangle & |4\rangle & |7\rangle & |6\rangle\\
  \hline
  \frac{1}{\sqrt2}(|1\rangle+|3\rangle) & \frac{1}{\sqrt2}(|0\rangle-|2\rangle) & \frac{1}{2}(|0\rangle+|1\rangle+|2\rangle-|3\rangle) & \frac{1}{2}(|0\rangle-|1\rangle+|2\rangle+|3\rangle) & \color{red!50!yellow}|6_6\rangle & \color{red!50!yellow}|7_6\rangle & \color{blue!30!green}|4_7\rangle & \color{blue!30!green}|5_7\rangle\\
  \hline
  \frac{1}{\sqrt2}(|1\rangle-|3\rangle) & \frac{1}{\sqrt2}(|0\rangle+|2\rangle) & \frac{1}{2}(|0\rangle+|1\rangle-|2\rangle+|3\rangle) & \frac{1}{2}(|0\rangle-|1\rangle-|2\rangle-|3\rangle) & \color{red!50!yellow}|7_6\rangle & \color{red!50!yellow}|6_6\rangle & \color{blue!30!green}|5_7\rangle & \color{blue!30!green}|4_7\rangle\\
  \arrayrulecolor{red}\hline
  \color{red}|4^0\rangle & \color{red}|5^0\rangle & \color{red}|6^0\rangle & \color{red}|7^0\rangle & \color{magenta}|0_4\rangle & \color{magenta}|1_4\rangle & \color{magenta}|2_4\rangle & |3\rangle\\
  \arrayrulecolor{black}\hline
  \color{red}|5^1\rangle & \color{red}|4^1\rangle & \color{red}|7^1\rangle & \color{red}|6^1\rangle & \color{magenta}|1_4\rangle  & \color{magenta}|0_4\rangle & |3\rangle & \color{magenta}|2_4\rangle\\
  \hline
  \color{red}|6^2\rangle & \color{red}|7^2\rangle & \color{red}|4^2\rangle & \color{red}|5^2\rangle & \color{magenta}|2_4\rangle & |3\rangle & \color{magenta}|0_4\rangle & \color{magenta}|1_4\rangle\\
  \hline
  \color{red}|7^3\rangle & \color{red}|6^3\rangle & \color{red}|5^3\rangle & \color{red}|4^3\rangle & |3\rangle & \color{magenta}|2_4\rangle & \color{magenta}|1_4\rangle & \color{magenta}|0_4\rangle\\
  \hline
\end{array}
$}
\end{center}

where
$\color{magenta}|0_4\rangle ,~\color{magenta}|1_4\rangle ,~\color{magenta}|2_4\rangle$ are given by Case 1,
$\color{red!50!yellow}|6_6\rangle ,~\color{red!50!yellow}|7_6\rangle$ are given by Case 3,
$\color{blue!30!green}|4_7\rangle ,~\color{blue!30!green}|5_7\rangle$ are given by Case 4 and
 $\color{red}|4^i\rangle,~\color{red}|5^i\rangle,~\color{red}|6^i\rangle,~\color{red}|7^i\rangle$  are given by Case 7 for $0\leq i \leq 3$.

\item There exists a QLS(8) with $c$ = 44.\label{44}

\begin{center}
\vspace{0.2cm}\centering\scalebox{0.8}{
$
\Psi_{44} =
\begin{array}{!{\color{black}\vline}c!{\color{black}\vline}c!{\color{black}\vline}c!{\color{black}\vline}c!{\color{red}\vline}c!{\color{black}\vline}c!{\color{black}\vline}c!{\color{black}\vline}c!{\color{black}\vline}}
  \hline
  |0\rangle & |1\rangle & \frac{1}{\sqrt2}(|2\rangle-|3\rangle) & \frac{1}{\sqrt2}(|2\rangle+|3\rangle) & |4\rangle & |5\rangle & |6\rangle & |7\rangle\\
  \hline
  |2\rangle & |3\rangle & \frac{1}{\sqrt2}(|0\rangle-|1\rangle) & \frac{1}{\sqrt2}(|0\rangle+|1\rangle) & |5\rangle & |4\rangle & |7\rangle & |6\rangle\\
  \hline
  \frac{1}{\sqrt2}(|1\rangle+|3\rangle) & \frac{1}{\sqrt2}(|0\rangle-|2\rangle) & \frac{1}{2}(|0\rangle+|1\rangle+|2\rangle-|3\rangle) & \frac{1}{2}(|0\rangle-|1\rangle+|2\rangle+|3\rangle) & \color{red!50!yellow}|6_6\rangle & \color{red!50!yellow}|7_6\rangle & \color{blue!30!green}|4_7\rangle & \color{blue!30!green}|5_7\rangle\\
  \hline
  \frac{1}{\sqrt2}(|1\rangle-|3\rangle) & \frac{1}{\sqrt2}(|0\rangle+|2\rangle) & \frac{1}{2}(|0\rangle+|1\rangle-|2\rangle+|3\rangle) & \frac{1}{2}(|0\rangle-|1\rangle-|2\rangle-|3\rangle) & \color{red!50!yellow}|7_6\rangle & \color{red!50!yellow}|6_6\rangle & \color{blue!30!green}|5_7\rangle & \color{blue!30!green}|4_7\rangle\\
  \arrayrulecolor{red}\hline
  \color{red}|4^0\rangle & \color{red}|5^0\rangle & \color{red}|6^0\rangle & \color{red}|7^0\rangle & \color{violet}|0_{10}\rangle & \color{violet}|1_{10}\rangle & \color{violet}|2_{10}\rangle & \color{violet}|3_{10}\rangle\\
  \arrayrulecolor{black}\hline
  \color{red}|5^1\rangle & \color{red}|4^1\rangle & \color{red}|7^1\rangle & \color{red}|6^1\rangle & \color{violet}|1_{10}\rangle  & \color{violet}|0_{10}\rangle & \color{violet}|3_{10}\rangle & \color{violet}|2_{10}\rangle\\
  \hline
  \color{red}|6^2\rangle & \color{red}|7^2\rangle & \color{red}|4^2\rangle & \color{red}|5^2\rangle & \color{violet}|2_{10}\rangle & \color{violet}|3_{10}\rangle & \color{violet}|0_{10}\rangle & \color{violet}|1_{10}\rangle\\
  \hline
  \color{red}|7^3\rangle & \color{red}|6^3\rangle & \color{red}|5^3\rangle & \color{red}|4^3\rangle & \color{violet}|3_{10}\rangle & \color{violet}|2_{10}\rangle & \color{violet}|1_{10}\rangle & \color{violet}|0_{10}\rangle\\
  \hline
\end{array}
$}
\end{center}

where
$\color{red!50!yellow}|6_6\rangle ,~\color{red!50!yellow}|7_6\rangle$ are given by Case 3,
$\color{blue!30!green}|4_7\rangle ,~\color{blue!30!green}|5_7\rangle$ are given by Case 4,
$\color{violet}|0_{10}\rangle ,~ \color{violet}|1_{10}\rangle ,~ \color{violet}|2_{10}\rangle ,~ \color{violet}|3_{10}\rangle$ are given by Case 10 and
 $\color{red}|4^i\rangle,~\color{red}|5^i\rangle,~\color{red}|6^i\rangle,~\color{red}|7^i\rangle$ are given by Case 7 for $0\leq i \leq 3$.

\item There exists a QLS(8) with $c$ = 46.\label{46}

\begin{center}
\vspace{0.2cm}\centering\scalebox{0.8}{
$
\Psi_{46} =
\begin{array}{!{\color{black}\vline}c!{\color{black}\vline}c!{\color{black}\vline}c!{\color{black}\vline}c!{\color{red}\vline}c!{\color{black}\vline}c!{\color{black}\vline}c!{\color{black}\vline}c!{\color{black}\vline}}
  \hline
  |0\rangle & |1\rangle & |2\rangle & |3\rangle & |4\rangle & |5\rangle & |6\rangle & |7\rangle\\
  \hline
  |1\rangle & |0\rangle & |3\rangle & |2\rangle & |5\rangle & |4\rangle & |7\rangle & |6\rangle\\
  \hline
  \color{red!60!black}|2_8\rangle & \color{red!60!black}|3_8\rangle & \color{yellow!50!black}|0_9\rangle & \color{yellow!50!black}|1_9\rangle & \color{red!50!yellow}|6_6\rangle & \color{red!50!yellow}|7_6\rangle & |4\rangle & |5\rangle\\
  \hline
  \color{red!60!black}|3_8\rangle & \color{red!60!black}|2_8\rangle & \color{yellow!50!black}|1_9\rangle & \color{yellow!50!black}|0_9\rangle & \color{red!50!yellow}|7_6\rangle & \color{red!50!yellow}|6_6\rangle & |5\rangle & |4\rangle\\
  \arrayrulecolor{red}\hline
  \color{red}|4^0\rangle & \color{red}|5^0\rangle & \color{red}|6^0\rangle & \color{red}|7^0\rangle & \color{violet!50!yellow}|0^0\rangle & \color{violet!50!yellow}|1^0\rangle & \color{violet!50!yellow}|2^0\rangle & \color{violet!50!yellow}|3^0\rangle\\
  \arrayrulecolor{black}\hline
  \color{red}|5^1\rangle & \color{red}|4^1\rangle & \color{red}|7^1\rangle & \color{red}|6^1\rangle & \color{violet!50!yellow}|1^1\rangle  & \color{violet!50!yellow}|0^1\rangle & \color{violet!50!yellow}|3^1\rangle & \color{violet!50!yellow}|2^1\rangle\\
  \hline
  \color{red}|6^2\rangle & \color{red}|7^2\rangle & \color{red}|4^2\rangle & \color{red}|5^2\rangle & \color{violet!50!yellow}|2^2\rangle & \color{violet!50!yellow}|3^2\rangle & \color{violet!50!yellow}|0^2\rangle & \color{violet!50!yellow}|1^2\rangle\\
  \hline
  \color{red}|7^3\rangle & \color{red}|6^3\rangle & \color{red}|5^3\rangle & \color{red}|4^3\rangle & \color{violet!50!yellow}|3^3\rangle & \color{violet!50!yellow}|2^3\rangle & \color{violet!50!yellow}|1^3\rangle & \color{violet!50!yellow}|0^3\rangle\\
  \hline
\end{array}
$}
\end{center}

where
$\color{red!50!yellow}|6_6\rangle ,~  \color{red!50!yellow}|7_6\rangle$ are given by Case 3,
$\color{red!60!black}|2_8\rangle ,~ \color{red!60!black}|3_8\rangle $ are given by Case 5,
$\color{yellow!50!black}|0_9\rangle ,~ \color{yellow!50!black}|1_9\rangle$ are given by Case 6,
$\color{violet}|0_{10}\rangle ,~ \color{violet}|1_{10}\rangle ,~ \color{violet}|2_{10}\rangle ,~ \color{violet}|3_{10}\rangle$ are given by Case 10 and
$\color{red}|4^i\rangle,~\color{red}|5^i\rangle,~\color{red}|6^i\rangle,~\color{red}|7^i\rangle$ are given by Case 7 for $0\leq i \leq 3$.
Moreover, let $\Phi = (|\Phi_{ij}\rangle)=(\alpha^i_j)$ be a QLS$(4)$ as defined in Case 7, where $i,j \in [4]$.
For $4\leq i,~j \leq 7$, define ${|\Psi_{46}\rangle}_{ij} = $
$
\begin{pmatrix}
\alpha^{i-4}_{j-4}\\
{{\mathbf 0}_{4}}
\end{pmatrix}
$, i.e.,
$$\begin{array}{lll}
\color{violet!50!yellow}|0^0\rangle \color{black} = \frac{1}{2}(|0\rangle + |1\rangle + |2\rangle + |3\rangle),  &
\color{violet!50!yellow}|1^0\rangle \color{black} =\frac{1}{2}(|0\rangle - \mathrm{i}|1\rangle - |2\rangle + \mathrm{i}|3\rangle),&
\color{violet!50!yellow}|2^0\rangle \color{black} = \frac{1}{2}(|0\rangle - |1\rangle + |2\rangle - |3\rangle),\\
\color{violet!50!yellow}|3^0\rangle \color{black} =\frac{1}{2}(|0\rangle + \mathrm{i}|1\rangle - |2\rangle - \mathrm{i}|3\rangle),&
\color{violet!50!yellow}|1^1\rangle \color{black} = \frac{1}{2}(|0\rangle - |1\rangle - \mathrm{i}|2\rangle + \mathrm{i}|3\rangle),&
\color{violet!50!yellow}|0^1\rangle \color{black} = \frac{1}{2}(|0\rangle + \mathrm{i}|1\rangle + \mathrm{i}|2\rangle - |3\rangle),\\
\color{violet!50!yellow}|3^1\rangle \color{black} = \frac{1}{2}(|0\rangle + |1\rangle - \mathrm{i}|2\rangle - \mathrm{i}|3\rangle),&
\color{violet!50!yellow}|2^1\rangle \color{black} = \frac{1}{2}(|0\rangle - \mathrm{i}|1\rangle + \mathrm{i}|2\rangle + |3\rangle),&
\color{violet!50!yellow}|2^2\rangle \color{black} = \frac{1}{2}(|0\rangle + |1\rangle - |2\rangle - |3\rangle),\\
\color{violet!50!yellow}|3^2\rangle \color{black} = \frac{1}{2}(|0\rangle - \mathrm{i}|1\rangle + |2\rangle - \mathrm{i}|3\rangle),&
\color{violet!50!yellow}|0^2\rangle \color{black} = \frac{1}{2}(|0\rangle - |1\rangle - |2\rangle + |3\rangle),&
\color{violet!50!yellow}|1^2\rangle \color{black} = \frac{1}{2}(|0\rangle + \mathrm{i}|1\rangle + |2\rangle + \mathrm{i}|3\rangle),\\
\color{violet!50!yellow}|3^3\rangle \color{black} = \frac{1}{2}(|0\rangle - |1\rangle + \mathrm{i}|2\rangle - \mathrm{i}|3\rangle),&
\color{violet!50!yellow}|2^3\rangle \color{black} = \frac{1}{2}(|0\rangle + \mathrm{i}|1\rangle - \mathrm{i}|2\rangle + |3\rangle),&
\color{violet!50!yellow}|1^3\rangle \color{black} = \frac{1}{2}(|0\rangle + |1\rangle + \mathrm{i}|2\rangle + \mathrm{i}|3\rangle),\\
\color{violet!50!yellow}|0^3\rangle \color{black} = \frac{1}{2}(|0\rangle - \mathrm{i}|1\rangle - \mathrm{i}|2\rangle - |3\rangle).&&
\end{array}$$

\item There exists a QLS(8) with $c$ = 47.\label{47}

\begin{center}
\vspace{0.2cm}\centering\scalebox{0.8}{
$
\Psi_{47} =
\begin{array}{!{\color{black}\vline}c!{\color{black}\vline}c!{\color{black}\vline}c!{\color{black}\vline}c!{\color{red}\vline}c!{\color{black}\vline}c!{\color{black}\vline}c!{\color{black}\vline}c!{\color{black}\vline}}
  \hline
  |0\rangle & |1\rangle & \frac{1}{\sqrt2}(|2\rangle-|3\rangle) & \frac{1}{\sqrt2}(|2\rangle+|3\rangle) & \color{teal}|4^0\rangle & \color{teal}|5^0\rangle & \color{teal}|6^0\rangle & \color{teal}|7^0\rangle\\
  \hline
  |2\rangle & |3\rangle & \frac{1}{\sqrt2}(|0\rangle-|1\rangle) & \frac{1}{\sqrt2}(|0\rangle+|1\rangle) & \color{teal}|5^1\rangle & \color{teal}|4^1\rangle & \color{teal}|7^1\rangle & \color{teal}|6^1\rangle\\
  \hline
  \frac{1}{\sqrt2}(|1\rangle+|3\rangle) & \frac{1}{\sqrt2}(|0\rangle-|2\rangle) & \frac{1}{2}(|0\rangle+|1\rangle+|2\rangle-|3\rangle) & \frac{1}{2}(|0\rangle-|1\rangle+|2\rangle+|3\rangle) & \color{teal}|6^2\rangle & \color{teal}|7^2\rangle & \color{teal}|4^2\rangle & \color{teal}|5^2\rangle\\
  \hline
  \frac{1}{\sqrt2}(|1\rangle-|3\rangle) & \frac{1}{\sqrt2}(|0\rangle+|2\rangle) & \frac{1}{2}(|0\rangle+|1\rangle-|2\rangle+|3\rangle) & \frac{1}{2}(|0\rangle-|1\rangle-|2\rangle-|3\rangle) & \color{teal}|7^3\rangle & \color{teal}|6^3\rangle & \color{teal}|5^3\rangle & \color{teal}|4^3\rangle\\
  \arrayrulecolor{red}\hline
  |4\rangle & |5\rangle & \frac{1}{\sqrt2}(|6\rangle-|7\rangle) & \frac{1}{\sqrt2}(|6\rangle+|7\rangle) & \color{magenta}|0_4\rangle & \color{magenta}|1_4\rangle & \color{magenta}|2_4\rangle & |3\rangle\\
  \arrayrulecolor{black}\hline
  |6\rangle & |7\rangle & \frac{1}{\sqrt2}(|0\rangle-|1\rangle) & \frac{1}{\sqrt2}(|4\rangle+|5\rangle) & \color{magenta}|1_4\rangle  & \color{magenta}|0_4\rangle & |3\rangle & \color{magenta}|2_4\rangle\\
  \hline
  \frac{1}{\sqrt2}(|5\rangle+|7\rangle) & \frac{1}{\sqrt2}(|4\rangle-|6\rangle) & \frac{1}{2}(|4\rangle+|5\rangle+|6\rangle-|7\rangle) & \frac{1}{2}(|4\rangle-|5\rangle+|6\rangle+|7\rangle) & \color{magenta}|2_4\rangle & |3\rangle & \color{magenta}|0_4\rangle & \color{magenta}|1_4\rangle\\
  \hline
  \frac{1}{\sqrt2}(|5\rangle-|7\rangle) & \frac{1}{\sqrt2}(|4\rangle+|6\rangle) & \frac{1}{2}(|4\rangle+|5\rangle-|6\rangle+|7\rangle) & \frac{1}{2}(|4\rangle-|5\rangle-|6\rangle-|7\rangle) & |3\rangle & \color{magenta}|2_4\rangle & \color{magenta}|1_4\rangle & \color{magenta}|0_4\rangle\\
  \hline
\end{array}
$}
\end{center}

where
$\color{magenta}|0_4\rangle ,~\color{magenta}|1_4\rangle ,~\color{magenta}|2_4\rangle$ are given by Case 1. And
let
$U'_i = F^{\dagger}_{4} \cdot U_i = (\beta^{i}_{0}~\beta^{i}_{1}~ \beta^{i}_{2} ~\beta^{i}_{3})$ with $U_i$ mentioned in Eq.(\ref{uu}) of order 4. For $0\leq i \leq 3,~4\leq j \leq 7$, define
${|\Psi_{47}\rangle}_{ij} = $
$
\begin{pmatrix}
{{\mathbf 0}_{4}}\\
\beta^{i}_{j-4}
\end{pmatrix}
$, i.e., for any $4 \leq i \leq 7$, $\color{teal}|i^0\rangle~\color{black}= ~|i\rangle$ and
$$\begin{array}{lll}
\color{teal}|4^1\rangle \color{black} = \frac{1}{2}[ |5\rangle + (1-\mathrm{i})|6\rangle + \mathrm{i}|7\rangle],&
\color{teal}|5^1\rangle \color{black} =\frac{1}{2}[\mathrm{i}|4\rangle  + |6\rangle + (1-\mathrm{i})|7\rangle],&
\color{teal}|6^1\rangle \color{black} = \frac{1}{2}[(1-\mathrm{i})|4\rangle - \mathrm{i}|5\rangle + |7\rangle],\\
\color{teal}|7^1\rangle \color{black} =\frac{1}{2}[|4\rangle + (1-\mathrm{i})|5\rangle - \mathrm{i}|6\rangle],&
\color{teal}|5^2\rangle \color{black} = \frac{1}{2}[(1+\mathrm{i})|5\rangle + (1-\mathrm{i})|7\rangle],&
\color{teal}|4^2\rangle \color{black} = \frac{1}{2}[(1-\mathrm{i})|4\rangle + (1+\mathrm{i})|6\rangle ],\\
\color{teal}|7^3\rangle \color{black} = \frac{1}{2}[(1-\mathrm{i})|5\rangle + (1+\mathrm{i})|7\rangle],&
\color{teal}|6^3\rangle \color{black} = \frac{1}{2}[(1+\mathrm{i})|4\rangle + (1-\mathrm{i})|6\rangle ],&
\color{teal}|6^3\rangle \color{black} = \frac{1}{2}[-\mathrm{i}|5\rangle + (1+\mathrm{i})|6\rangle - |7\rangle],\\
\color{teal}|7^3\rangle \color{black} = \frac{1}{2}[|4\rangle -\mathrm{i}|6\rangle + (1+\mathrm{i})|7\rangle],&
\color{teal}|4^2\rangle \color{black} = \frac{1}{2}[(1+\mathrm{i})|4\rangle + |5\rangle -\mathrm{i}|7\rangle],&
\color{teal}|5^2\rangle \color{black} = \frac{1}{2}[-\mathrm{i}|4\rangle + (1+\mathrm{i})|5\rangle + |6\rangle].
\end{array}$$

\item There exists a QLS(8) with $c$ = 48.\label{48}

\begin{center}
\vspace{0.2cm}\centering\scalebox{0.8}{
$
\Psi_{48} =
\begin{array}{!{\color{black}\vline}c!{\color{black}\vline}c!{\color{black}\vline}c!{\color{black}\vline}c!{\color{red}\vline}c!{\color{black}\vline}c!{\color{black}\vline}c!{\color{black}\vline}c!{\color{black}\vline}}
  \hline
  |0\rangle & |1\rangle & |2\rangle & |3\rangle & |4\rangle & |5\rangle & |6\rangle & |7\rangle\\
  \hline
  |1\rangle & |0\rangle & |3\rangle & |2\rangle & |5\rangle & |4\rangle & |7\rangle & |6\rangle\\
  \hline
  \color{red!60!black}|2_8\rangle & \color{red!60!black}|3_8\rangle & \color{yellow!50!black}|0_9\rangle & \color{yellow!50!black}|1_9\rangle & \color{red!50!yellow}|6_6\rangle & \color{red!50!yellow}|7_6\rangle & \color{blue!30!green}|4_7\rangle & \color{blue!30!green}|5_7\rangle\\
  \hline
  \color{red!60!black}|3_8\rangle & \color{red!60!black}|2_8\rangle & \color{yellow!50!black}|1_9\rangle & \color{yellow!50!black}|0_9\rangle & \color{red!50!yellow}|7_6\rangle & \color{red!50!yellow}|6_6\rangle & \color{blue!30!green}|5_7\rangle & \color{blue!30!green}|4_7\rangle\\
  \arrayrulecolor{red}\hline
  \color{red}|4^0\rangle & \color{red}|5^0\rangle & \color{red}|6^0\rangle & \color{red}|7^0\rangle & \color{violet!50!yellow}|0^0\rangle & \color{violet!50!yellow}|1^0\rangle & \color{violet!50!yellow}|2^0\rangle & \color{violet!50!yellow}|3^0\rangle\\
  \arrayrulecolor{black}\hline
  \color{red}|5^1\rangle & \color{red}|4^1\rangle & \color{red}|7^1\rangle & \color{red}|6^1\rangle & \color{violet!50!yellow}|1^1\rangle  & \color{violet!50!yellow}|0^1\rangle & \color{violet!50!yellow}|3^1\rangle & \color{violet!50!yellow}|2^1\rangle\\
  \hline
  \color{red}|6^2\rangle & \color{red}|7^2\rangle & \color{red}|4^2\rangle & \color{red}|5^2\rangle & \color{violet!50!yellow}|2^2\rangle & \color{violet!50!yellow}|3^2\rangle & \color{violet!50!yellow}|0^2\rangle & \color{violet!50!yellow}|1^2\rangle\\
  \hline
  \color{red}|7^3\rangle & \color{red}|6^3\rangle & \color{red}|5^3\rangle & \color{red}|4^3\rangle & \color{violet!50!yellow}|3^3\rangle & \color{violet!50!yellow}|2^3\rangle & \color{violet!50!yellow}|1^3\rangle & \color{violet!50!yellow}|0^3\rangle\\
  \hline
\end{array}
$}
\end{center}

where
$\color{red!50!yellow}|6_6\rangle ,~  \color{red!50!yellow}|7_6\rangle$ are given by Case 3,
$\color{blue!30!green}|4_7\rangle ,~  \color{blue!30!green}|5_7\rangle$ are given by Case 4,
$\color{red!60!black}|2_8\rangle ,~ \color{red!60!black}|3_8\rangle $ are given by Case 5,
$\color{yellow!50!black}|0_9\rangle ,~ \color{yellow!50!black}|1_9\rangle$ are given by Case 6,
$\color{red}|4^i\rangle,~\color{red}|5^i\rangle,~\color{red}|6^i\rangle,~\color{red}|7^i\rangle$ are given by Case 7
and $\color{violet!50!yellow}|0^i\rangle,~\color{violet!50!yellow}|1^i\rangle,
~\color{violet!50!yellow}|2^i\rangle,~\color{violet!50!yellow}|3^i\rangle$ are given by Case 21  for $0\leq i \leq 3$.

\item There exists a QLS(8) with $c$ = 50.\label{50}

\begin{center}
\vspace{0.2cm}\centering\scalebox{0.8}{
$
\Psi_{50} =
\begin{array}{!{\color{black}\vline}c!{\color{black}\vline}c!{\color{black}\vline}c!{\color{black}\vline}c!{\color{red}\vline}c!{\color{black}\vline}c!{\color{black}\vline}c!{\color{black}\vline}c!{\color{black}\vline}}
  \hline
  |0\rangle & |1\rangle & \frac{1}{\sqrt2}(|2\rangle-|3\rangle) & \frac{1}{\sqrt2}(|2\rangle+|3\rangle) & |4\rangle & |5\rangle & |6\rangle & |7\rangle\\
  \hline
  |2\rangle & |3\rangle & \frac{1}{\sqrt2}(|0\rangle-|1\rangle) & \frac{1}{\sqrt2}(|0\rangle+|1\rangle) & |5\rangle & |4\rangle & |7\rangle & |6\rangle\\
  \hline
  \frac{1}{\sqrt2}(|1\rangle+|3\rangle) & \frac{1}{\sqrt2}(|0\rangle-|2\rangle) & \frac{1}{2}(|0\rangle+|1\rangle+|2\rangle-|3\rangle) & \frac{1}{2}(|0\rangle-|1\rangle+|2\rangle+|3\rangle) & \color{red!50!yellow}|6_6\rangle & \color{red!50!yellow}|7_6\rangle & |4\rangle & |5\rangle\\
  \hline
  \frac{1}{\sqrt2}(|1\rangle-|3\rangle) & \frac{1}{\sqrt2}(|0\rangle+|2\rangle) & \frac{1}{2}(|0\rangle+|1\rangle-|2\rangle+|3\rangle) & \frac{1}{2}(|0\rangle-|1\rangle-|2\rangle-|3\rangle) & \color{red!50!yellow}|7_6\rangle & \color{red!50!yellow}|6_6\rangle & |5\rangle & |4\rangle\\
  \arrayrulecolor{red}\hline
  |4\rangle & |5\rangle & \frac{1}{\sqrt2}(|6\rangle-|7\rangle) & \frac{1}{\sqrt2}(|6\rangle+|7\rangle) & \color{violet!50!yellow}|0^0\rangle & \color{violet!50!yellow}|1^0\rangle & \color{violet!50!yellow}|2^0\rangle & \color{violet!50!yellow}|3^0\rangle\\
  \arrayrulecolor{black}\hline
  |6\rangle & |7\rangle & \frac{1}{\sqrt2}(|0\rangle-|1\rangle) & \frac{1}{\sqrt2}(|4\rangle+|5\rangle) & \color{violet!50!yellow}|1^1\rangle  & \color{violet!50!yellow}|0^1\rangle & \color{violet!50!yellow}|3^1\rangle & \color{violet!50!yellow}|2^1\rangle\\
  \hline
  \frac{1}{\sqrt2}(|5\rangle+|7\rangle) & \frac{1}{\sqrt2}(|4\rangle-|6\rangle) & \frac{1}{2}(|4\rangle+|5\rangle+|6\rangle-|7\rangle) & \frac{1}{2}(|4\rangle-|5\rangle+|6\rangle+|7\rangle) & \color{violet!50!yellow}|2^2\rangle & \color{violet!50!yellow}|3^2\rangle & \color{violet!50!yellow}|0^2\rangle & \color{violet!50!yellow}|1^2\rangle\\
  \hline
  \frac{1}{\sqrt2}(|5\rangle-|7\rangle) & \frac{1}{\sqrt2}(|4\rangle+|6\rangle) & \frac{1}{2}(|4\rangle+|5\rangle-|6\rangle+|7\rangle) & \frac{1}{2}(|4\rangle-|5\rangle-|6\rangle-|7\rangle) & \color{violet!50!yellow}|3^3\rangle & \color{violet!50!yellow}|2^3\rangle & \color{violet!50!yellow}|1^3\rangle & \color{violet!50!yellow}|0^3\rangle\\
  \hline
\end{array}
$}
\end{center}

where
$\color{red!50!yellow}|6_6\rangle ,~  \color{red!50!yellow}|7_6\rangle$ are given by Case 3 and $\color{violet!50!yellow}|0^i\rangle,~\color{violet!50!yellow}|1^i\rangle,
~\color{violet!50!yellow}|2^i\rangle,~\color{violet!50!yellow}|3^i\rangle$ are given by Case 21 for $0\leq i \leq 3$.

\item There exists a QLS(8) with $c$ = 51.\label{51}

\begin{center}
\vspace{0.2cm}\centering\scalebox{0.8}{
$
\Psi_{51} =
\begin{array}{!{\color{black}\vline}c!{\color{black}\vline}c!{\color{black}\vline}c!{\color{black}\vline}c!{\color{red}\vline}c!{\color{black}\vline}c!{\color{black}\vline}c!{\color{black}\vline}c!{\color{black}\vline}}
  \hline
  |0\rangle & |1\rangle & \frac{1}{\sqrt2}(|2\rangle-|3\rangle) & \frac{1}{\sqrt2}(|2\rangle+|3\rangle) & \color{red}|4^0\rangle & \color{red}|5^0\rangle & \color{red}|6^0\rangle & \color{red}|7^0\rangle\\
  \hline
  |2\rangle & |3\rangle & \frac{1}{\sqrt2}(|0\rangle-|1\rangle) & \frac{1}{\sqrt2}(|0\rangle+|1\rangle) & \color{red}|5^1\rangle & \color{red}|4^1\rangle & \color{red}|7^1\rangle & \color{red}|6^1\rangle\\
  \hline
  \frac{1}{\sqrt2}(|1\rangle+|3\rangle) & \frac{1}{\sqrt2}(|0\rangle-|2\rangle) & \frac{1}{2}(|0\rangle+|1\rangle+|2\rangle-|3\rangle) & \frac{1}{2}(|0\rangle-|1\rangle+|2\rangle+|3\rangle) & \color{red}|6^2\rangle & \color{red}|7^2\rangle & \color{red}|4^2\rangle & \color{red}|5^2\rangle\\
  \hline
  \frac{1}{\sqrt2}(|1\rangle-|3\rangle) & \frac{1}{\sqrt2}(|0\rangle+|2\rangle) & \frac{1}{2}(|0\rangle+|1\rangle-|2\rangle+|3\rangle) & \frac{1}{2}(|0\rangle-|1\rangle-|2\rangle-|3\rangle) & \color{red}|7^3\rangle & \color{red}|6^3\rangle & \color{red}|5^3\rangle & \color{red}|4^3\rangle\\
  \arrayrulecolor{red}\hline
  |4\rangle & |5\rangle & \frac{1}{\sqrt2}(|6\rangle-|7\rangle) & \frac{1}{\sqrt2}(|6\rangle+|7\rangle) & \color{magenta}|0_4\rangle & \color{magenta}|1_4\rangle & \color{magenta}|2_4\rangle & |3\rangle\\
  \arrayrulecolor{black}\hline
  |6\rangle & |7\rangle & \frac{1}{\sqrt2}(|0\rangle-|1\rangle) & \frac{1}{\sqrt2}(|4\rangle+|5\rangle) & \color{magenta}|1_4\rangle  & \color{magenta}|0_4\rangle & |3\rangle & \color{magenta}|2_4\rangle\\
  \hline
  \frac{1}{\sqrt2}(|5\rangle+|7\rangle) & \frac{1}{\sqrt2}(|4\rangle-|6\rangle) & \frac{1}{2}(|4\rangle+|5\rangle+|6\rangle-|7\rangle) & \frac{1}{2}(|4\rangle-|5\rangle+|6\rangle+|7\rangle) & \color{magenta}|2_4\rangle & |3\rangle & \color{magenta}|0_4\rangle & \color{magenta}|1_4\rangle\\
  \hline
  \frac{1}{\sqrt2}(|5\rangle-|7\rangle) & \frac{1}{\sqrt2}(|4\rangle+|6\rangle) & \frac{1}{2}(|4\rangle+|5\rangle-|6\rangle+|7\rangle) & \frac{1}{2}(|4\rangle-|5\rangle-|6\rangle-|7\rangle) & |3\rangle & \color{magenta}|2_4\rangle & \color{magenta}|1_4\rangle & \color{magenta}|0_4\rangle\\
  \hline
\end{array}
$}
\end{center}

where
$\color{magenta}|0_4\rangle ,~\color{magenta}|1_4\rangle ,~\color{magenta}|2_4\rangle$ are given by Case 1 and
$\color{red}|4^i\rangle,~\color{red}|5^i\rangle,~\color{red}|6^i\rangle,~\color{red}|7^i\rangle$ are given by Case 7 for $0\leq i \leq 3$.

\item There exists a QLS(8) with $c$ = 52.\label{52}

\begin{center}
\vspace{0.2cm}\centering\scalebox{0.8}{
$
\Psi_{52} =
\begin{array}{!{\color{black}\vline}c!{\color{black}\vline}c!{\color{black}\vline}c!{\color{black}\vline}c!{\color{red}\vline}c!{\color{black}\vline}c!{\color{black}\vline}c!{\color{black}\vline}c!{\color{black}\vline}}
  \hline
  |0\rangle & |1\rangle & \frac{1}{\sqrt2}(|2\rangle-|3\rangle) & \frac{1}{\sqrt2}(|2\rangle+|3\rangle) & |4\rangle & |5\rangle & |6\rangle & |7\rangle\\
  \hline
  |2\rangle & |3\rangle & \frac{1}{\sqrt2}(|0\rangle-|1\rangle) & \frac{1}{\sqrt2}(|0\rangle+|1\rangle) & |5\rangle & |4\rangle & |7\rangle & |6\rangle\\
  \hline
  \frac{1}{\sqrt2}(|1\rangle+|3\rangle) & \frac{1}{\sqrt2}(|0\rangle-|2\rangle) & \frac{1}{2}(|0\rangle+|1\rangle+|2\rangle-|3\rangle) & \frac{1}{2}(|0\rangle-|1\rangle+|2\rangle+|3\rangle) & \color{red!50!yellow}|6_6\rangle & \color{red!50!yellow}|7_6\rangle & \color{blue!30!green}|4_7\rangle & \color{blue!30!green}|5_7\rangle\\
  \hline
  \frac{1}{\sqrt2}(|1\rangle-|3\rangle) & \frac{1}{\sqrt2}(|0\rangle+|2\rangle) & \frac{1}{2}(|0\rangle+|1\rangle-|2\rangle+|3\rangle) & \frac{1}{2}(|0\rangle-|1\rangle-|2\rangle-|3\rangle) & \color{red!50!yellow}|7_6\rangle & \color{red!50!yellow}|6_6\rangle & \color{blue!30!green}|5_7\rangle & \color{blue!30!green}|4_7\rangle\\
  \arrayrulecolor{red}\hline
  |4\rangle & |5\rangle & \frac{1}{\sqrt2}(|6\rangle-|7\rangle) & \frac{1}{\sqrt2}(|6\rangle+|7\rangle) & \color{violet!50!yellow}|0^0\rangle & \color{violet!50!yellow}|1^0\rangle & \color{violet!50!yellow}|2^0\rangle & \color{violet!50!yellow}|3^0\rangle\\
  \arrayrulecolor{black}\hline
  |6\rangle & |7\rangle & \frac{1}{\sqrt2}(|0\rangle-|1\rangle) & \frac{1}{\sqrt2}(|4\rangle+|5\rangle) & \color{violet!50!yellow}|1^1\rangle  & \color{violet!50!yellow}|0^1\rangle & \color{violet!50!yellow}|3^1\rangle & \color{violet!50!yellow}|2^1\rangle\\
  \hline
  \frac{1}{\sqrt2}(|5\rangle+|7\rangle) & \frac{1}{\sqrt2}(|4\rangle-|6\rangle) & \frac{1}{2}(|4\rangle+|5\rangle+|6\rangle-|7\rangle) & \frac{1}{2}(|4\rangle-|5\rangle+|6\rangle+|7\rangle) & \color{violet!50!yellow}|2^2\rangle & \color{violet!50!yellow}|3^2\rangle & \color{violet!50!yellow}|0^2\rangle & \color{violet!50!yellow}|1^2\rangle\\
  \hline
  \frac{1}{\sqrt2}(|5\rangle-|7\rangle) & \frac{1}{\sqrt2}(|4\rangle+|6\rangle) & \frac{1}{2}(|4\rangle+|5\rangle-|6\rangle+|7\rangle) & \frac{1}{2}(|4\rangle-|5\rangle-|6\rangle-|7\rangle) & \color{violet!50!yellow}|3^3\rangle & \color{violet!50!yellow}|2^3\rangle & \color{violet!50!yellow}|1^3\rangle & \color{violet!50!yellow}|0^3\rangle\\
  \hline
\end{array}
$}
\end{center}

where
$\color{red!50!yellow}|6_6\rangle ,~  \color{red!50!yellow}|7_6\rangle$ are given by Case 3,
$\color{blue!30!green}|4_7\rangle ,~  \color{blue!30!green}|5_7\rangle$ are given by Case 4 and
$\color{violet!50!yellow}|0^i\rangle,~\color{violet!50!yellow}|1^i\rangle,
~\color{violet!50!yellow}|2^i\rangle,~\color{violet!50!yellow}|3^i\rangle$ are given by Case 21 for $0\leq i \leq 3$.

\item There exists a QLS(8) with $c$ = 54.\label{54}

\begin{center}
\vspace{0.2cm}\centering\scalebox{0.8}{
$
\Psi_{54} =
\begin{array}{!{\color{black}\vline}c!{\color{black}\vline}c!{\color{black}\vline}c!{\color{black}\vline}c!{\color{red}\vline}c!{\color{black}\vline}c!{\color{black}\vline}c!{\color{black}\vline}c!{\color{black}\vline}}
  \hline
  |0\rangle & |1\rangle & \frac{1}{\sqrt2}(|2\rangle-|3\rangle) & \frac{1}{\sqrt2}(|2\rangle+|3\rangle) & |4\rangle & |5\rangle & |6\rangle & |7\rangle\\
  \hline
  |2\rangle & |3\rangle & \frac{1}{\sqrt2}(|0\rangle-|1\rangle) & \frac{1}{\sqrt2}(|0\rangle+|1\rangle) & |5\rangle & |4\rangle & |7\rangle & |6\rangle\\
  \hline
  \frac{1}{\sqrt2}(|1\rangle+|3\rangle) & \frac{1}{\sqrt2}(|0\rangle-|2\rangle) & \frac{1}{2}(|0\rangle+|1\rangle+|2\rangle-|3\rangle) & \frac{1}{2}(|0\rangle-|1\rangle+|2\rangle+|3\rangle) & \color{red!50!yellow}|6_6\rangle & \color{red!50!yellow}|7_6\rangle & |4\rangle & |5\rangle\\
  \hline
  \frac{1}{\sqrt2}(|1\rangle-|3\rangle) & \frac{1}{\sqrt2}(|0\rangle+|2\rangle) & \frac{1}{2}(|0\rangle+|1\rangle-|2\rangle+|3\rangle) & \frac{1}{2}(|0\rangle-|1\rangle-|2\rangle-|3\rangle) & \color{red!50!yellow}|7_6\rangle & \color{red!50!yellow}|6_6\rangle & |5\rangle & |4\rangle\\
  \arrayrulecolor{red}\hline
  \color{red}|4^0\rangle & \color{red}|5^0\rangle & \color{red}|6^0\rangle & \color{red}|7^0\rangle & \color{violet!50!yellow}|0^0\rangle & \color{violet!50!yellow}|1^0\rangle & \color{violet!50!yellow}|2^0\rangle & \color{violet!50!yellow}|3^0\rangle\\
  \arrayrulecolor{black}\hline
  \color{red}|5^1\rangle & \color{red}|4^1\rangle & \color{red}|7^1\rangle & \color{red}|6^1\rangle & \color{violet!50!yellow}|1^1\rangle  & \color{violet!50!yellow}|0^1\rangle & \color{violet!50!yellow}|3^1\rangle & \color{violet!50!yellow}|2^1\rangle\\
  \hline
  \color{red}|6^2\rangle & \color{red}|7^2\rangle & \color{red}|4^2\rangle & \color{red}|5^2\rangle & \color{violet!50!yellow}|2^2\rangle & \color{violet!50!yellow}|3^2\rangle & \color{violet!50!yellow}|0^2\rangle & \color{violet!50!yellow}|1^2\rangle\\
  \hline
  \color{red}|7^3\rangle & \color{red}|6^3\rangle & \color{red}|5^3\rangle & \color{red}|4^3\rangle & \color{violet!50!yellow}|3^3\rangle & \color{violet!50!yellow}|2^3\rangle & \color{violet!50!yellow}|1^3\rangle & \color{violet!50!yellow}|0^3\rangle\\
  \hline
\end{array}
$}
\end{center}

where
$\color{red!50!yellow}|6_6\rangle ,~  \color{red!50!yellow}|7_6\rangle$ are given by Case 3,
$\color{red}|4^i\rangle,~\color{red}|5^i\rangle,~\color{red}|6^i\rangle,~\color{red}|7^i\rangle$ are given by Case 7 and $\color{violet!50!yellow}|0^i\rangle,~\color{violet!50!yellow}|1^i\rangle,
~\color{violet!50!yellow}|2^i\rangle,~\color{violet!50!yellow}|3^i\rangle$ are given by Case 21 for $0\leq i \leq 3$..

\item There exists a QLS(8) with $c$ = 56.\label{56}

\begin{center}
\vspace{0.2cm}\centering\scalebox{0.8}{
$
\Psi_{56} =
\begin{array}{!{\color{black}\vline}c!{\color{black}\vline}c!{\color{black}\vline}c!{\color{black}\vline}c!{\color{red}\vline}c!{\color{black}\vline}c!{\color{black}\vline}c!{\color{black}\vline}c!{\color{black}\vline}}
  \hline
  |0\rangle & |1\rangle & \frac{1}{\sqrt2}(|2\rangle-|3\rangle) & \frac{1}{\sqrt2}(|2\rangle+|3\rangle) & |4\rangle & |5\rangle & |6\rangle & |7\rangle\\
  \hline
  |2\rangle & |3\rangle & \frac{1}{\sqrt2}(|0\rangle-|1\rangle) & \frac{1}{\sqrt2}(|0\rangle+|1\rangle) & |5\rangle & |4\rangle & |7\rangle & |6\rangle\\
  \hline
  \frac{1}{\sqrt2}(|1\rangle+|3\rangle) & \frac{1}{\sqrt2}(|0\rangle-|2\rangle) & \frac{1}{2}(|0\rangle+|1\rangle+|2\rangle-|3\rangle) & \frac{1}{2}(|0\rangle-|1\rangle+|2\rangle+|3\rangle) & \color{red!50!yellow}|6_6\rangle & \color{red!50!yellow}|7_6\rangle & \color{blue!30!green}|4_7\rangle & \color{blue!30!green}|5_7\rangle\\
  \hline
  \frac{1}{\sqrt2}(|1\rangle-|3\rangle) & \frac{1}{\sqrt2}(|0\rangle+|2\rangle) & \frac{1}{2}(|0\rangle+|1\rangle-|2\rangle+|3\rangle) & \frac{1}{2}(|0\rangle-|1\rangle-|2\rangle-|3\rangle) & \color{red!50!yellow}|7_6\rangle & \color{red!50!yellow}|6_6\rangle & \color{blue!30!green}|5_7\rangle & \color{blue!30!green}|4_7\rangle\\
  \arrayrulecolor{red}\hline
  \color{red}|4^0\rangle & \color{red}|5^0\rangle & \color{red}|6^0\rangle & \color{red}|7^0\rangle & \color{violet!50!yellow}|0^0\rangle & \color{violet!50!yellow}|1^0\rangle & \color{violet!50!yellow}|2^0\rangle & \color{violet!50!yellow}|3^0\rangle\\
  \arrayrulecolor{black}\hline
  \color{red}|5^1\rangle & \color{red}|4^1\rangle & \color{red}|7^1\rangle & \color{red}|6^1\rangle & \color{violet!50!yellow}|1^1\rangle  & \color{violet!50!yellow}|0^1\rangle & \color{violet!50!yellow}|3^1\rangle & \color{violet!50!yellow}|2^1\rangle\\
  \hline
  \color{red}|6^2\rangle & \color{red}|7^2\rangle & \color{red}|4^2\rangle & \color{red}|5^2\rangle & \color{violet!50!yellow}|2^2\rangle & \color{violet!50!yellow}|3^2\rangle & \color{violet!50!yellow}|0^2\rangle & \color{violet!50!yellow}|1^2\rangle\\
  \hline
  \color{red}|7^3\rangle & \color{red}|6^3\rangle & \color{red}|5^3\rangle & \color{red}|4^3\rangle & \color{violet!50!yellow}|3^3\rangle & \color{violet!50!yellow}|2^3\rangle & \color{violet!50!yellow}|1^3\rangle & \color{violet!50!yellow}|0^3\rangle\\
  \hline
\end{array}
$}
\end{center}

where
$\color{red!50!yellow}|6_6\rangle ,~  \color{red!50!yellow}|7_6\rangle$ are given by Case 3,
$\color{blue!30!green}|4_7\rangle ,~  \color{blue!30!green}|5_7\rangle$ are given by Case 4,
$\color{red}|4^i\rangle,~\color{red}|5^i\rangle,~\color{red}|6^i\rangle,~\color{red}|7^i\rangle$ are given by Case 7 and
$\color{violet!50!yellow}|0^i\rangle,~\color{violet!50!yellow}|1^i\rangle,
~\color{violet!50!yellow}|2^i\rangle,~\color{violet!50!yellow}|3^i\rangle$ are given by Case 21 for $0\leq i \leq 3$.

\item There exists a QLS(8) with $c$ = 60.\label{60}

\begin{center}
\vspace{0.2cm}\centering\scalebox{0.8}{
$
\Psi_{60} =
\begin{array}{!{\color{black}\vline}c!{\color{black}\vline}c!{\color{black}\vline}c!{\color{black}\vline}c!{\color{red}\vline}c!{\color{black}\vline}c!{\color{black}\vline}c!{\color{black}\vline}c!{\color{black}\vline}}
  \hline
  |0\rangle & |1\rangle & \frac{1}{\sqrt2}(|2\rangle-|3\rangle) & \frac{1}{\sqrt2}(|2\rangle+|3\rangle) & \color{teal}|4^0\rangle & \color{teal}|5^0\rangle & \color{teal}|6^0\rangle & \color{teal}|7^0\rangle\\
  \hline
  |2\rangle & |3\rangle & \frac{1}{\sqrt2}(|0\rangle-|1\rangle) & \frac{1}{\sqrt2}(|0\rangle+|1\rangle) & \color{teal}|5^1\rangle & \color{teal}|4^1\rangle & \color{teal}|7^1\rangle & \color{teal}|6^1\rangle\\
  \hline
  \frac{1}{\sqrt2}(|1\rangle+|3\rangle) & \frac{1}{\sqrt2}(|0\rangle-|2\rangle) & \frac{1}{2}(|0\rangle+|1\rangle+|2\rangle-|3\rangle) & \frac{1}{2}(|0\rangle-|1\rangle+|2\rangle+|3\rangle) & \color{teal}|6^2\rangle & \color{teal}|7^2\rangle & \color{teal}|4^2\rangle & \color{teal}|5^2\rangle\\
  \hline
  \frac{1}{\sqrt2}(|1\rangle-|3\rangle) & \frac{1}{\sqrt2}(|0\rangle+|2\rangle) & \frac{1}{2}(|0\rangle+|1\rangle-|2\rangle+|3\rangle) & \frac{1}{2}(|0\rangle-|1\rangle-|2\rangle-|3\rangle) & \color{teal}|7^3\rangle & \color{teal}|6^3\rangle & \color{teal}|5^3\rangle & \color{teal}|4^3\rangle\\
  \arrayrulecolor{red}\hline
  |4\rangle & |5\rangle & \frac{1}{\sqrt2}(|6\rangle-|7\rangle) & \frac{1}{\sqrt2}(|6\rangle+|7\rangle) & \color{violet!50!yellow}|0^0\rangle & \color{violet!50!yellow}|1^0\rangle & \color{violet!50!yellow}|2^0\rangle & \color{violet!50!yellow}|3^0\rangle\\
  \arrayrulecolor{black}\hline
  |6\rangle & |7\rangle & \frac{1}{\sqrt2}(|0\rangle-|1\rangle) & \frac{1}{\sqrt2}(|4\rangle+|5\rangle) & \color{violet!50!yellow}|1^1\rangle  & \color{violet!50!yellow}|0^1\rangle & \color{violet!50!yellow}|3^1\rangle & \color{violet!50!yellow}|2^1\rangle\\
  \hline
  \frac{1}{\sqrt2}(|5\rangle+|7\rangle) & \frac{1}{\sqrt2}(|4\rangle-|6\rangle) & \frac{1}{2}(|4\rangle+|5\rangle+|6\rangle-|7\rangle) & \frac{1}{2}(|4\rangle-|5\rangle+|6\rangle+|7\rangle) & \color{violet!50!yellow}|2^2\rangle & \color{violet!50!yellow}|3^2\rangle & \color{violet!50!yellow}|0^2\rangle & \color{violet!50!yellow}|1^2\rangle\\
  \hline
  \frac{1}{\sqrt2}(|5\rangle-|7\rangle) & \frac{1}{\sqrt2}(|4\rangle+|6\rangle) & \frac{1}{2}(|4\rangle+|5\rangle-|6\rangle+|7\rangle) & \frac{1}{2}(|4\rangle-|5\rangle-|6\rangle-|7\rangle) & \color{violet!50!yellow}|3^3\rangle & \color{violet!50!yellow}|2^3\rangle & \color{violet!50!yellow}|1^3\rangle & \color{violet!50!yellow}|0^3\rangle\\
  \hline
\end{array}
$}
\end{center}

where $\color{violet!50!yellow}|0^i\rangle,~\color{violet!50!yellow}|1^i\rangle,
~\color{violet!50!yellow}|2^i\rangle,~\color{violet!50!yellow}|3^i\rangle$ are given by Case 21 and
$\color{teal}|4^i\rangle,~\color{teal}|5^i\rangle,~\color{teal}|6^i\rangle,~\color{teal}|7^i\rangle$ are given by Case 22 for $0\leq i \leq 3$.

\end{enumerate}
\end{document}